\documentclass[11pt]{amsart}

\usepackage{preamble-locally-gentle}

\let\amp=&
\usepackage{makecell}

\begin{document}

\title{Homological conditions on locally gentle algebras}

\address{Department of Mathematics, Box 354350, University of Washington, Seattle, Washington 98195, USA}

\author{S. Ford}
\email{sarafinaf@unr.edu}
\author{A. Oswald}
\email{amreio@uw.edu}
\author{J. J. Zhang}
\email{zhang@math.washington.edu}

\begin{abstract}
Gentle algebras are a class of special biserial algebras whose representation theory has been thoroughly described. In this paper, we consider the infinite dimensional generalizations of gentle algebras, referred to as locally gentle algebras. We give combinatorial descriptions of the center, prime spectrum, and homological dimensions of a locally gentle algebra, including an explicit injective resolution. We classify when these algebras are Artin-Schelter Gorenstein,  Artin-Schelter regular, and Cohen-Macaulay, and provide an analogue of Stanley's theorem for locally gentle algebras.
\end{abstract}



\maketitle

\section{Introduction}\label{sec:intro}
Gentle algebras are a well-understood class of finite dimensional bound quiver algebras that were originally studied as iterated tilted algebras of type $A$ \cite{AH82} and $\tilde{A}$ \cite{AsSk87}. Their indecomposable module categories are described in \cite{BSi21,HKK17,OPS18}. They are closed under derived equivalence \cite{SZ03}.
The indecomposable objects of their derived categories are classified in \cite{BM03}, and \cite{ALP} gave a basis for the space of morphisms between those indecomposable complexes. The cones of those morphisms were described in \cite{CPS19}, and their almost-split triangles were described in \cite{Bo11}.  They are known to be Koszul \cite{BH08} and Gorenstein \cite{GR05}, and their global dimension and self-injective dimension can be calculated combinatorially \cite{LGH}. 

Locally gentle algebras are an infinite dimensional generalization of gentle algebras (see Definition \ref{def:gentle}), which arise naturally as the Koszul duals of gentle algebras \cite{BH08}. Recently, there has been increased interest in this class of algebra. Their Hochschild cohomology is described in \cite{CSSSA24} and \cite{PPP19} gives a bijection between the isomorphism classes and homeomorphism classes of certain marked surfaces. The torsion classes in the category of finite length modules over locally gentle algebras are classified in \cite{CD20}. Finitely generated modules of string algebras, which include locally gentle algebras, are described in \cite{Crawley-Boevey18} in terms of string and band modules.

In this paper, we determine when locally gentle algebras satisfy certain homological conditions. Along the way, we extend known results about the homological dimensions of gentle algebras. 
In Section \ref{sec:prelim}, we give the definitions and conventions used throughout the paper, and detail some preliminary properties of locally gentle algebras. One key invariant of a gentle algebra, $A=\kk\cQ/\cI$, is that each arrow in $\cQ_1$ is contained in a unique maximal path, i.e., a path $w$ in $\cQ$ so that $bw=wb=0$ for any arrow $b\in\cQ_1$.
The notion of maximal paths carries over to locally gentle algebras in the sense that  
each arrow is contained in a unique path which is either maximal in the sense of gentle algebras (referred to as a \emph{finite maximal path}) or an infinite maximal path which  is determined by a primitive cycle $c$ which is unique up to rotation and whose powers are nonzero. An immediate consequence of this is Proposition \ref{prop:h_A}, which states that the Hilbert series of a (locally) gentle algebra $A=\kk\cQ/\cI$ is
\begin{equation*}
    h_A(t)=\frac{\abs{\cQ_0}+(\abs{\cQ_1}-\abs{\cQ_0})t-\sum_{n\geq 1}k_nt^{n+1}}{1-t},
\end{equation*}
where $\cQ$ is a finite quiver and $k_n$ is the number of finite maximal paths of length at least $n$. Additionally, the center of a (locally) gentle algebra is generated by certain finite maximal paths and elements $m_\gamma$ corresponding to infinite maximal paths (see Lemma \ref{lem:center}).

In Section \ref{sec:primes}, we give a description of the prime spectrum of a (locally) gentle algebra in terms of infinite maximal paths and the central elements $m_\gamma$ and conclude that a (locally) gentle algebra is prime if and only if it has a unique, infinite maximal path.

\begin{thm}[Theorem \ref{thm:primeideals}] Let $A$ be a (locally) gentle algebra over a base field $\kk$. The prime ideals of $A$ are exactly the following. 
\begin{enumerate}
    \item For each $v\in\cQ_0$, the maximal ideal $\fm_v$ of $A$ which does not contain $v$.
    \item For each infinite maximal path $\gamma$, the annihilator of the ideal $\cJ_\gamma$ generated by nonstationary subpaths of $\gamma$.
    \item  For each infinite maximal path $\gamma$ and an irreducible polynomial $p(t)\in \kk[t]$ with $\deg(p)\geq 1$ and $p(0)\neq 0$, the ideal 
    \begin{equation*}\mathfrak{p}(\gamma,p)=\anni_A(\cJ_\gamma)+\langle p(m_\gamma)\rangle\end{equation*}
\end{enumerate}
	
If $\gamma$ is an infinite maximal path in $A$ and $\kk$ is algebraically closed, we have the following inclusions.
\begin{center}
    \begin{tikzpicture}
        \node at (0,-.25) {$\anni_A(\cJ_\gamma)$};
        \node at (-2,1.25) {$\{\fm_v:\gamma \text{  passes through }v\}$};
        \node at (2.5,1.25) {$\{\mathfrak{p}(\gamma,t-\lambda)\}_{\lambda\in\kk^\times}$};
        \draw[-] (-2,1) to node[midway, below, rotate=-35]{$\supseteq$} (-.75,.2);
        \draw[-] (2,1) to node[midway, below, rotate=35]{$\subseteq$} (.75,.2);
    \end{tikzpicture}
\end{center}
\end{thm}

Homological dimensions of gentle and almost gentle algebras have been described in \cite{LGH, WHL25}. In particular, \cite{LGH} gave a description of the global and self-injective dimensions of a gentle algebra in terms of paths in its Koszul dual.  
In Section \ref{sec:homological-dim}, we extend these descriptions to locally gentle algebras. Let $A^\#$ be the Koszul dual of $A$ and let $\cL'$ be the set of finite maximal paths of $A^\#$.
\begin{thm}[Theorems \ref{thm:gldim} and \ref{thm:injdim}]
Let $A=\kk\cQ/\cI$ be a (locally) gentle algebra. Then 
\begin{equation*}
    \gldim(A)=\max_{\alpha\in\cQ_1}\ell(\gamma^\#_\alpha),
\end{equation*}
where $\gamma^\#_\alpha$ is the maximal path in $A^\#$ containing the arrow $\alpha$ and
\begin{equation*}
    \injdim(A_A)=\injdim(_AA)=\begin{cases}
    \max_{p\in \cL'}\ell(p) & \text{if $\cL'\neq \emptyset$},\\
    0 & \text{if $A=\kk\widetilde{A}_n/(\kk\widetilde{A}_n)_{\geq 2}$},\\
    1 & \text{otherwise}.
    \end{cases}
\end{equation*}
\end{thm}
In particular, like gentle algebras, locally gentle algebras are always Iwanaga--Gorenstein. 
In both \cite{LGH} and \cite{GR05}, which first established that gentle algebras are Gorenstein, the results about injective dimension used the projective dimension and relied on the correspondence between injective modules over $A$ and projective modules over $A^{\op}$. In particular, $\injdim_A(A)=\pdim_{A^{\op}}(DA)$, which is not necessarily true in the infinite dimensional case. Instead, we provide explicit injective resolutions in Lemma \ref{lem:injres}.

In \cite{RR19}, it was shown that a locally finite graded elementary $\kk$-algebra (in particular, any (locally) gentle algebra) is a twisted Calabi-Yau algebra of dimension $d$ if and only if it is Artin-Schelter regular with global dimension $d$. Moreover, it was shown that $\kk\widetilde{A}_n$ is the only (locally) gentle twisted Calabi-Yau algebra of dimension 1 and hence the only (locally) gentle  Artin-Schelter regular algebra with global dimension 1. In Section \ref{sec:ASconds}, we show that $\kk\widetilde{A}_n$ is the only (locally) gentle Artin-Schelter regular algebra of any dimension. Additionally, in Proposition \ref{prop:locgentleAScond}, we show that a (locally) gentle algebra is  Artin-Schelter Gorenstein if the corresponding quiver is either $\widetilde{A}_n$ (with $A=\kk\widetilde{A}_n$ or $(\kk\widetilde{A}_n)^\#$) or has twice as many arrows as vertices.

In Section \ref{subsec:CM}, we consider two definitions of Cohen-Macaulay for a Noetherian $\mathbb{N}$-graded algebra which is not necessarily connected and find that both are equivalent for a (locally) gentle algebra.
\begin{thm}[Theorem \ref{thm:CM}]
For $A$ (locally) gentle, the following are equivalent:
\begin{itemize}
    \item  $\GKdim(A)=\dpth(A)$.
    \item  $A$ is a finitely generated, free $\kk[x_1,\ldots, x_d]$-module for some $x_i\in A$ and $d=\GKdim(A)$.
    \item  $A$ has no finite maximal paths or $A$ has no infinite maximal paths.
\end{itemize}
{When these equivalent conditions hold, we say that $A$ is {Cohen-Macaulay}.}
\end{thm}
This definition of Cohen-Macaulay allows an analogue of Stanley's theorem for locally gentle algebras, although this result fails for gentle algebras (see Example \ref{ex: Stanley fails}).
\begin{thm}[Theorem \ref{thm:Stanleylocgentlealgs}]
Let $A$ be a locally gentle algebra. If  $A$ is Cohen-Macaulay, then the following are equivalent:
\begin{itemize}
    \item[(i)]  $A$ is Artin-Schelter Gorenstein.
    \item[(ii)] $h_A(t^{-1})=\pm t^{k}h_A(t)$ for some  $k\in\ZZ$.
    \item[(iii)] $h_A(t^{-1})=\pm t^{k}h_A(t)$ for some  $k=0$ or $1$.
\end{itemize}
\end{thm}

\section{Preliminaries on locally gentle algebras}\label{sec:prelim}

Throughout this paper, we work over the base field $\kk$ with characteristic 0, which is not assumed to be algebraically closed unless otherwise stated. For a set $S$, we denote its $\kk$-linear span by $\kk S$.

\subsection{Graded algebras and modules}

In this subsection, we outline notation and definitions used throughout the paper, following \cite{toolkit} and \cite{RR19}. Let $\gA$ be a locally finite, left bounded graded algebra over $\kk$ which is not necessarily connected.  Then $\gA_i$ denotes the degree $i$ portion of $\gA$, so that 
$\gA=\bigoplus_{n\geq 0} \gA_n$. We let $\gA_{\geq n}:=\bigoplus_{i\geq n}\gA_i$ and similarly $\gA_{\leq n}:=\bigoplus_{i\leq n}\gA_i$. In particular, we often write $A_+$ to mean $A_{\geq 1}$. The \emph{Hilbert series of $\gA$} is the power series
\begin{equation*}
    h_\gA(t)=\sum_{n=0}^\infty \dim_{\kk}(\gA_n)t^n.
\end{equation*}

Let $V$ be a finite dimensional subspace which generates $\gA$ and let $V^n$ be the $\kk$-subspace of $\gA$ spanned by elements $v_1\cdots v_m$, where $v_i\in V$ and $m\leq n$. The \emph{Gelfand-Kirillov dimension}  of $\gA$ is
\begin{equation*}
\GKdim(\gA)=\limsup_{n\rightarrow\infty}\frac{\log(\dim_\kk(V^n))}{\log(n)}
\end{equation*} 
 (see \cite{GKdimbook}). This definition is independent of the choice of $V$. 

Let $\grmod{\gA}$ be the category of $\ZZ$-graded right $\gA$-modules. Throughout this paper, we work in this category and, unless otherwise specified, all $\gA$-modules are assumed to be $\ZZ$-graded right modules. In particular, when we refer to projective and injective $\gA$-modules, we mean projective and injective objects in $\grmod{\gA}$.
Let $M$ be a $\gA$-module. We write $M_i$ to denote the degree $i$ portion of $M$, so that $M=\bigoplus_{i\in\ZZ}M_i$. The \emph{$k$-th shift} of $M$, denoted $M[k]$, is the $\gA$-module with $(M[k])_i=M_{i+k}$.

For $\gA$-modules $M,N$, 
\begin{equation*}\grHom_\gA(M,N)\defeq \bigoplus \Hom_\gA^i(M,N),\end{equation*}
where $\Hom_\gA^i(M,N)$ are the graded homomorphisms of degree $i$, and $\grExt_\gA^j \defeq  R^j\grHom_\gA$.  
The \emph{graded Jacobson radical} of $\gA$, denoted $\grJ(\gA)$, is the intersection of its maximal homogeneous ideals. Then, we have that 
\begin{equation*}
\grJ(\gA)=J(\gA_0)\oplus \gA_+=J(\gA_0)\oplus \bigoplus_{i\geq 1}\gA_i,
\end{equation*}
where $J(\gA_0)$ is the Jacobson radical of $\gA_0$. We let $S\defeq \gA/\grJ(\gA)$.  The \emph{depth} of $\gA$, denoted $\dpth(\gA)$, is
\begin{equation*}
    \dpth(\gA)=\min\setst{i}{ \grExt_\gA^i(S,\gA)\neq 0}.
\end{equation*}

\subsection{Locally gentle algebras}

For a finite set $X$, let $\abs{X}$ denote the number of elements in $X$.
A \emph{finite quiver}, $\cQ=(\cQ_0,\cQ_1,s,t)$, is a directed graph where $\cQ_0$ is a set of vertices, $\cQ_1$ is a set of arrows, and the functions $s,t:\cQ_1\to\cQ_0$ give the \emph{source} and \emph{target} of each arrow respectively. Throughout, we assume that $1\leq \abs{\cQ_1},\abs{\cQ_0}<\infty$, and 
that any quiver is connected, meaning that the underlying graph is connected.

For a vertex $v\in\cQ_0$, let $\indeg(v)$ be the number of arrows $\alpha\in\cQ_1$ so that $t(\alpha)=v$, and let $\outdeg(v)$ be the number of arrows $\beta\in\cQ_1$ so that $s(\beta)=v$.
If $\indeg(v)=0$, then we say that $v$ is a \emph{source}, and if $\outdeg(v)=0$, then we say that $v$ is a \emph{sink}.

A \emph{path}, $p$, in a quiver is a concatenation of arrows $p=\alpha_0\alpha_1\dotsc\alpha_n$ where $t(\alpha_i)=s(\alpha_{i+1})$ for every $1\leq i <n$. We say the path $p$ has \emph{length} $n$, and we denote this by $\ell(p)=n$. We denote the first arrow of $p$ by $F(p)$, i.e. $F(p)=\alpha_1$, and the last arrow of $p$ by $L(p)$, i.e. $L(p)=\alpha_n$. Further, we let $s(p)=s(F(p))$ and $t(p)=t(L(p))$. Note, when we say $s^{-1}(v)$ and $t^{-1}(v)$ for some $v\in\cQ_0$, we mean all arrows with source $v$ and target $v$ respectively. 

For every $v\in\cQ_0$, we let $\epsilon_v$ be a stationary path at the vertex $v$, which has length 0 and $s(\epsilon_v)=t(\epsilon_v)\defeq v$. We let $F(\epsilon_v)=L(\epsilon_v)\defeq \epsilon_v$.
If $p$ and $q$ are two paths so that $t(p)=s(q)$, then $pq$ is the path formed by concatenating these two paths. If we have paths $r,q',q''$ so that $p=q'rq''$, then we say that $r$ is a \emph{subpath} of $p$ and denote this $r\leq p$. We say that $r$ is an \emph{initial} (\emph{terminal}) subpath of $p$ if $q'$ ($q''$) is stationary. If $X$ is a set of paths in $\cQ$, denote the subset of $X$ consisting of paths of length $n$ by $X_n$ and the subset of $X$ consisting of paths of positive length by $X_+$.

Given a quiver $\cQ$, we can form its \emph{path algebra} $\kk\cQ$, whose basis consists of all paths in $\cQ$ and whose multiplication is given by concatenation of paths. 
Any path algebra $\kk\cQ$ is graded by path length. We say that an ideal of $\kk\cQ$ is homogeneous if it is generated by homogeneous elements.
If $\cI$ is a homogeneous ideal of $\kk\cQ$, then the algebra $A=\kk\cQ/\cI$ is also graded by path length. If $\cQ$ has $n$ vertices, then $A_0\cong (\kk)^n$ as an algebra, so $J(A_0)$ is the zero ideal and $\grJ(A)=A_+$.
If the ideal $\cI$ is generated by paths, then the algebra $\kk\cQ/\cI$ is a monomial algebra with a basis consisting of paths. 

\begin{defn}\label{def:gentle} An algebra $A=\kk\cQ/\cI$ is called a \emph{gentle algebra} if
\begin{enumerate}[(i)]
    \item the quiver $\cQ$ is finite, connected, and $\indeg(v), \outdeg(v)\leq 2$ for each vertex $v\in\cQ_0$.
    \item whenever $s(\alpha)=t(\beta)=t(\beta')$ with $\beta\neq\beta'$, exactly one of $\beta\alpha$ and $\beta'\alpha$ is in $\cI$.
    \item whenever $t(\alpha)=s(\beta)=s(\beta')$ with $\beta\neq \beta'$, then exactly one of  $\alpha\beta$ and $\alpha\beta'$ is in $\cI$.
    \item the ideal $\cI$ is generated by paths of lengths 2.
    \item the algebra $A$ is finite dimensional (i.e. $\cI$ is admissible).
\end{enumerate}
We say the algebra $A$ is \emph{locally gentle} if it is infinite dimensional and satisfies conditions (i)-(iv). We write (locally) gentle to mean gentle or locally gentle. 
\end{defn}

Any (locally) gentle algebra is Noetherian \cite[Lemma 3.3]{Crawley-Boevey18}. However, a (locally) gentle algebra $A=\kk\cQ/\cI$ is Artinian if and only if it is gentle (as $A_{\geq n}$ is a descending chain of ideals). As the ideal $\cI$ is homogeneous and generated by paths, $A$ is a monomial algebra graded by path length and has a basis of paths. Note that the Gelfand-Kirillov dimension of a locally gentle algebra is 1.

\begin{notation}\label{not:algebra}
For the remainder of this paper, $\cQ$ denotes a finite, connected quiver, $\cI$ denotes an ideal in $\kk\cQ$ so that $\kk\cQ/\cI$ is (locally) gentle, $A$ denotes the (locally) gentle algebra $\kk\cQ/\cI$, and $\cB$ denotes the basis of $A$ consisting of paths.
\end{notation}

An important consequence of (ii) in the definition above is that for any path $p \in \cQ$, there is at most one arrow $\alpha$ so that $\alpha p\not\in\cI$. Similarly, there is at most one arrow $\beta\in\cQ_1$ so that $p\beta\not\in\cI$. It follows that if $F(w)=F(w')$ ($L(w)=L(w')$) for $w,w'\in\cB$ then either $w'$ is an initial (terminal) subpath of $w$ or vice versa. Additionally, since $\cI$ is generated by paths, $ww'=ww''\notin\cI$ implies $w'=w''$ for any paths $w,w',w''\in\cB$. 

Quivers where all vertices have the same in-degree and out-degree are of particular interest in this paper.

\begin{example}\label{ex:maxbiserial}
Suppose that $A=\kk\cQ/\cI$ is locally gentle. Then $\abs{\cQ_1}\leq 2\abs{\cQ_0}$ and $\abs{\cQ_1}=2\abs{\cQ_0}$ if and only if every vertex $v\in \cQ_0$ has $\indeg(v)=\outdeg(v)=2$ (i.e., the quiver is 2-regular). 
Such quivers with at most three vertices are shown below.
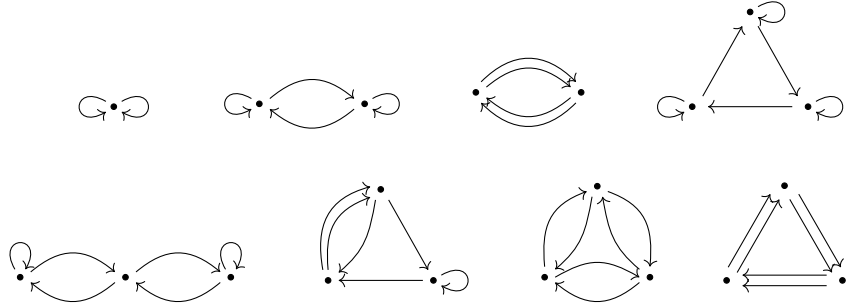
\begin{align*}
    \begin{tikzpicture}[scale=.7]
        \vtx{0,0}{}
        \draw[->] (.15,.1) to[out=30, in=-30, looseness=10] (.15,-.1) ;
        \draw[->] (-.15,.1) to[out=150, in=-150, looseness=10] (-.15,-.1) ;
    \end{tikzpicture} 
    \qquad
    \begin{tikzpicture}[scale=.7]
        \vtx{0,0}{}
        \vtx{2,0}{}
        \draw[->] (1.8,-.1) to[out=220, in=-40,looseness=1.2] (.2,-.1) ;
        \draw[->] (.2,.1) to[out=40, in=140,looseness=1.2] (1.8,.1) ;
        \draw[->] (-.15,.1) to[out=150, in=-150, looseness=10] (-.15,-.1) ;
        \draw[->] (2.15,.1) to[out=30, in=-30, looseness=10] (2.15,-.1) ;
    \end{tikzpicture}
    \qquad
    \begin{tikzpicture}[scale=.7]
        \vtx{0,0}{}
        \vtx{2,0}{}
        \draw[->] (1.9,-.2) to[out=225, in=-45,looseness=1.2] (.1,-.2) ;
        \draw[->] (1.8,-.1) to[out=220, in=-40,looseness=1.2]  (.2,-.1) ; 
        \draw[->] (.1,.2) to[out=45, in=135,looseness=1.2] (1.9,.2) ;
        \draw[->] (.2,.1) to[out=40, in=140,looseness=1.2] (1.8,.1) ; 
    \end{tikzpicture}
    \qquad
    \begin{tikzpicture}[scale=.7]
        \vtx{-.1,0}{}
        \vrx{1,1.8}{}
        \vtx{2.1,0}{}
        \draw[->] (.1,.2) to  (.85, 1.53);
        \draw[->] (1.15,1.53) to (1.9,.2);
        \draw[->] (1.8,0) to (.2,0);
        \draw[->] (-.25,.1) to[out=150, in=-150, looseness=10] (-.25,-.1) ;
        \draw[->] (1.15,1.9) to[out=30, in=-30, looseness=10]  (1.15,1.7) ;
        \draw[->] (2.25,.1) to[out=30, in=-30, looseness=10]  (2.25,-.1) ;
    \end{tikzpicture}
    \\
    \begin{tikzpicture}[scale=.7]
        \vtx{0,0}{}
        \vtx{2,0}{}
        \vtx{4,0}{}
        \draw[->] (3.8,-.1) to[out=220, in=-40,looseness=1.2] (2.2,-.1) ;
        \draw[->] (2.2,.1) to[out=40, in=140,looseness=1.2] (3.8,.1) ;
        \draw[->] (1.8,-.1) to[out=220, in=-40,looseness=1.2] (.2,-.1) ;
        \draw[->] (.2,.1) to[out=40, in=140,looseness=1.2] (1.8,.1) ;    
        \draw[->] (-.1,.15) to[out=120, in=60, looseness=10] (.1,.15) ;
        \draw[->] (3.9,.15) to[out=120, in=60, looseness=10]  (4.1,.15) ;
    \end{tikzpicture}
    \qquad
    \begin{tikzpicture}[scale=.7]
        \vtx{0,0}{}
        \vrx{1,1.73}{}
        \vtx{2,0}{}
        \draw[->] (.9,1.53) to[out=260, in=40,looseness=1.1] (.2,.15) ;
        \draw[->] (0,.2) to[out=95, in=200,looseness=1.2] (.8,1.58) ;
        \draw[->] (-.1,.2) to[out=100, in=195,looseness=1.2] (.8,1.78);
        \draw[->] (1.8,0) to (.2,0);
        \draw[->] (1.15,1.53) to (1.9,.2);
        \draw[->] (2.15,.1) to[out=30, in=-30, looseness=10] (2.15,-.1) ;
    \end{tikzpicture}
    \qquad
    \begin{tikzpicture}[scale=.7]
        \vtx{0,0}{}
        \vrx{1,1.73}{}
        \vtx{2,0}{}
        \draw[->] (0,.2) to[out=95, in=200,looseness=1.2] (.8,1.63) ;
        \draw[->] (.9,1.53) to[out=260, in=40,looseness=1.1] (.2,.15) ;
        \draw[->] (1.2,1.63) to[out=-20, in=85,looseness=1.2]  (2,.2);
        \draw[->] (1.8,.15)to[out=140, in=280,looseness=1.1] (1.1,1.53)  ;
        \draw[->] (1.8,-.1) to[out=220, in=-40,looseness=1.2] (.2,-.1) ;
        \draw[->] (.2,0) to[out=30, in=150,looseness=1.2]  (1.8,0) ;
    \end{tikzpicture}
    \qquad
    \begin{tikzpicture}[scale=.7]
        \vtx{-.1,0}{}
        \vrx{1,1.8}{}
        \vtx{2.1,0}{}
        \draw[->] (1.8,-.1) to  (.2,-.1);
        \draw[->] (1.8,.1) to  (.2,.1);
        \draw[->] (.1-.17, .15+.1) to   (.1-.17+.8,.15+.1+1.39);
        \draw[->] (.1,.15) to (.1+.8, .15+1.39);
        \draw[->] ({2-(.1-.17+.8)},.15+.1+1.39) to ({2-(.1-.17)}, .15+.1);
        \draw[->] ({2-(.1+.8)}, .15+1.39) to ({2-(.1)},.15);
    \end{tikzpicture}
\end{align*}
\end{example}

\begin{example}\label{ex:Antilde}
Suppose $\cQ$ is a quiver so that $\indeg(v)=\outdeg(v)=1$ for every $v\in\cQ_0$. Then, $\cQ=\widetilde{A}_{|\cQ_0|-1}$, where $\widetilde{A}_{n}$ is the extended Dynkin quiver of type $A$ with $n+1$ vertices. We denote the arrows of $\widetilde{A}_{n}$ by $\alpha_0,\alpha_1,\dotsc,\alpha_{n}$ where $s(\alpha_i)=i$ and $t(\alpha_i)=i+1 \,(\text{mod } n+1)$. Any choice of relations on $\widetilde{A}_n$ yields a (locally) gentle algebra. The quivers $\widetilde{A}_n$ for $n=0,1,2$ are shown below.
    
\begin{equation*}
    \begin{tikzpicture}[scale=.7]
        \vtx{0,0}{0}
        \draw[->] (-.1,.15) to[out=120, in=60, looseness=10] node[midway, above]{$\alpha_0$} (.1,.15);
    \end{tikzpicture} 
    \qquad
    \begin{tikzpicture}[scale=.7]
        \vtx{0,0}{0}
        \vtx{2,0}{1}
        \draw[->] (1.8,-.1) to[out=220, in=-40,looseness=1.2] node[midway, below]{$\alpha_1$} (.2,-.1) ;
        \draw[->] (.2,.1) to[out=40, in=140,looseness=1.2]node[midway, above]{$\alpha_0$} (1.8,.1) ;
    \end{tikzpicture}
    \qquad
    \begin{tikzpicture}[scale=.7]
        \vtx{-.1,0}{0}
        \vrx{1,1.8}{1}
        \vtx{2.1,0}{2}
        \draw[->] (.1,.2) to node[midway, left]{$\alpha_0$}  (.85, 1.53);
        \draw[->] (1.15,1.53) to node[midway, right]{$\alpha_1$} (1.9,.2);
        \draw[->] (1.8,0) to node[midway, below]{$\alpha_2$} (.2,0);
    \end{tikzpicture}
    \captionsetup{labelformat=empty}
\end{equation*}
\end{example}

\subsection{Maximal paths}\label{subsec:maxpaths}

Maximal paths will play an important role throughout the rest of this paper. In this section, we define the notion of an infinite maximal path and note some important properties.

\begin{defn}\label{def:maxpaths} 
Let $p$ be a path in $\cB$. We say that
\begin{itemize}
\item $p$ is a \emph{right maximal} if for every $\alpha\in\cQ_1$, $p\alpha\in\cI$.
\item $p$ is \emph{left maximal} if for every $\alpha\in\cQ_1$, $\alpha p\in\cI$.
\item $p$ is a \emph{finite maximal path} if $p$ is both left and right maximal. 
\end{itemize}
\end{defn}
\begin{example}\label{ex:nofinitemax}
     If $\cQ$ is 2-regular, then $A$ has no right maximal paths: for any path $p\in\cB$ there are two arrows $\alpha_1,\alpha_2$ with source $t(p)$ and only one of $p\alpha_1,p\alpha_2$ belongs to $\cI$. In particular, any (locally) gentle algebra on $\cQ$ has no finite maximal paths.
\end{example}

The following lemma is a consequence of conditions (ii) and (iii) in Definition \ref{def:gentle}.
\begin{lemma}
    For an arrow $\alpha\in A$, let $\gamma_r(\alpha)$ ($\gamma_l(\alpha)$) denote the unique path of maximal length in $\cB$ with no repeated arrows and first (last) arrow $\alpha$. Then, exactly one of the following holds.
    \begin{enumerate}
    \item The paths $\gamma_r(\alpha)$ and $\gamma_l(\alpha)$ are right and left maximal, respectively.
    \item The paths $\gamma_r(\alpha)$ and $\gamma_l(\alpha)$ are rotations of the same cycle and $\left(\gamma_r(\alpha)\right)^2, \left(\gamma_l(\alpha)\right)^2\notin\cI$.
    \end{enumerate}
\end{lemma}
\begin{proof}
    Suppose that $\gamma_r(\alpha)$ is not right maximal. Then, there exists an arrow $\beta$ such that $\gamma_r(\alpha)\beta\notin\cI$. By definition of $\gamma_r(\alpha)$, it must be that $\beta$ occurs in $\gamma_r(\alpha)$. That is, $\gamma_r(\alpha)=\alpha_1\cdots \alpha_n$ for some arrows $\alpha_i$ and $\beta=\alpha_j$ for some $j$. Now, since $\alpha_n\beta\notin\cI$, it follows that that either $j=1$ or $\alpha_{j-1}=\alpha_n$. The latter contradictions our definition of $\gamma_r(\alpha)$, so $\beta=\alpha_1=\alpha$. Hence $\alpha_n\alpha\notin\cI$. It follows that $\gamma_l(\alpha)=\alpha_2\cdots \alpha_n\alpha_1$ and $\left(\gamma_r(\alpha)\right)^2,\left(\gamma_l(\alpha)\right)^2\notin\cI$.
\end{proof}

In the first case, the arrow $\alpha$ (and any path containing $\alpha$) is contained in a unique finite maximal path, of which $\gamma_r(\alpha)$ and $\gamma_l(\alpha)$ are terminal and initial subpaths, respectively. Note that, as a consequence, no path in $\cB$ can contain two or more instances of $\alpha$.

In the second case, any path containing $\alpha$ must be a subpath of $\left(\gamma_r(\alpha)\right)^k$ for some $k>0$. This motivates the following definition.

\begin{defn}
A cycle $c$ in $A$ is \emph{primitive} if there is no cycle $c'$ in $A$ so that $c=(c')^k$ for $k>1$. Given a primitive cycle $c$ such that $c^2\not\in\cI$, we say that $c^\infty$ is an \emph{infinite maximal path} of $A$ and $c$ is a \emph{corresponding primitive cycle}. 
If $c'$ is a rotation of the cycle $c$, then $c'$ is also primitive and $(c')^2\notin\cI$ and we say the corresponding infinite maximal paths are equal ($(c')^\infty=c^\infty$). 

For a path $p$, we say that $p$ is a subpath of $c^\infty$ and write $p\leq c^\infty$ if $p\leq c^k$ for some $k>0$.
\end{defn}

\begin{corollary}\label{cor:maxpaths}
Fix an arrow $\alpha\in\cQ_1$. Then $\alpha$ is contained in a unique maximal path, which we denote by $\gamma_\alpha$. If $\gamma_r(\alpha)$ is right maximal, then $\gamma_\alpha$ is finite. Otherwise, $\gamma_\alpha=\left(\gamma_r(\alpha)\right)^\infty$.
\end{corollary}

If there is any path $p$ in $\cB_+$ containing arrows $\alpha$ and $\beta$, then $\gamma_\alpha=\gamma_\beta$ and $p\leq \gamma_\alpha$. Thus, the set $\cB_+$ is partitioned into subpaths of maximal paths. 
We note the following property of locally gentle algebras without finite maximal paths for  reference in future sections.
 
\begin{lemma}\label{lem:invoutdeg}
If $\cB$ contains no finite maximal paths, we have that $\indeg(v)=\outdeg(v)$ for all $v\in\cQ_0$, and exactly one of the following holds.
\begin{enumerate}
    \item The quiver is $\widetilde{A}_n$ and $A=\kk\widetilde{A}_n$ ($\indeg(v)=\outdeg(v)=1$ for each vertex $v$).
    \item The quiver is 2-regular ($\indeg(v)=\outdeg(v)=2$ for each vertex $v$).
    \item There is an arrow $\alpha$ with source $v$ and target $w$ such that $\indeg(v)=1$ and $\indeg(w)=2$.
\end{enumerate}
\end{lemma}

\begin{defn}
For a (locally) gentle algebra $A=\kk\cQ/\cI$, its Koszul dual is the (locally) gentle algebra $A^\#\defeq \kk\cQ/\cI^\#$ where, for any path $\alpha\beta$ of length two in $\cQ$, we have $\alpha\beta\in\cI^\#$ if and only if $\alpha\beta\notin \cI$ \cite{BH08}. 
We denote the basis of paths of $A^\#$ by $\cB^\#$, and let  $\gamma_\alpha^\#$ denote the (infinite) maximal path of $A^\#$ containing the arrow $\alpha\in\cQ_1$. Similarly, $\gamma_r^\#(\alpha)$ denotes the right maximal path or primitive cycle in $\cB^\#$ with $F(\gamma_r^\#(\alpha))=\alpha$, and $\gamma_\ell^\#(\alpha)$ denotes the left maximal path or primitive cycle in $\cB^\#$ with $L(\gamma_\ell^\#(\alpha))=\alpha$.
\end{defn}

For example, $\kk[x,y]/\langle xy\rangle$ and $\kk\langle x,y\rangle/ \langle x^2,y^2\rangle$ are Koszul dual, as are $\kk\widetilde{A}_n$ and $ \kk\widetilde{A}_n/(\kk\widetilde{A}_n)_{\geq 2}.$

\begin{defn}
The paths in $A^\#$ (i.e., elements of $\cB^\#$) are referred to as \emph{forbidden paths} of $A$. We say that a \emph{forbidden path} is  (right, left) maximal if it is  (right, left) maximal as a path in $A^\#$.
\end{defn}

\subsection{The Hilbert series}
In this section, we will describe the Hilbert series and show that 
\begin{equation*}
    \pm t^k h_A(t)=h_A(t^{-1})
\end{equation*}
holds only under certain conditions on the quiver $\cQ$, a result which we will use in Section \ref{sec:ASconds} for our analogue of Stanley's theorem (Theorem \ref{thm:Stanleylocgentlealgs}).

\begin{proposition}\label{prop:h_A}
Suppose that $A$ has $k_n$ finite maximal paths of length at least $n$. Then 
\begin{align*}
    h_A(t)& = \frac{\abs{\cQ_0}+(\abs{\cQ_1}-\abs{\cQ_0})t-\sum_{n\geq 1}k_nt^{n+1}}{1-t}.
\end{align*}
\end{proposition}

\begin{proof} Any path in $\cB_+$ is a subpath of some maximal path. An infinite maximal path containing $\ell$ distinct arrows has exactly $\ell$ subpaths of length $n$ for $1\leq n$, while a finite maximal path of length $\ell$  has  $\ell+1-n$ subpaths of length $n$ for $1\leq n\leq \ell$. 

Suppose that $A$ has $m$ finite maximal paths of lengths $\ell_1\leq\dotsc \leq \ell_m$. 
Since each arrow appears in at most one finite maximal path and finite maximal paths have no repeated arrows, there are $k\defeq \abs{\cQ_1}-\sum_{i=1}^m\ell_i$ arrows not contained in the $m$ finite maximal paths. Thus 
\begin{align*}
    h_A(t)& = 
    \abs{\cQ_0}+k\sum_{n=1}^\infty t^n+\sum_{i=1}^m\sum_{n=1}^{\ell_i}(\ell_i+1-n)t^n =\abs{\cQ_0}-k+\frac{k}{1-t}+\sum_{i=1}^m\sum_{n=1}^{\ell_i}(\ell_i+1-n)t^n.
\end{align*}
Observe that $(1-t)\sum_{n=1}^{\ell_i}(\ell_i+1-n)t^n=\ell_it-t\sum_{n=1}^{\ell_i}t^n$. Then $k_n=\max\setst{i}{\ell_i\geq n}$ for $1\leq n\leq \ell_m$ is the coefficient of $t^{n+1}$ in $\sum_{i=1}^m\sum_{n=1}^{\ell_i}t^{n+1}$. The result follows.
\end{proof}
\begin{example}\label{ex:hilb}
Let $A=\kk\cQ/\langle \alpha_1c, c\alpha_2, \beta_1\beta_2\rangle$, where $\cQ$ is the quiver
\begin{center}
    \begin{tikzpicture}[scale=.7]
    \vtx{0,0}{1}
    \vtx{2,0}{2}
    \vtx{4,0}{3}
    \vtx{6,0}{4}
    \vtx{8,0}{5}
    \draw[<-] (.2,0) to node[midway, above]{$\alpha_2$}(1.8,0) ;
    \draw[<-] (2.2,0) to node[midway, above]{$\alpha_1$}(3.8,0) ;
    \draw[<-] (1.9,.15) to[out=120, in=60, looseness=10] node[midway, above]{$c$}(2.1,.15) ;
    \draw[<-] (5.8,0) to node[midway, above]{$\beta_1$}(4.2,0) ;
    \draw[<-] (7.8,0) to node[midway, above]{$\beta_2$}(6.2,0) ;
    \end{tikzpicture}.
\end{center}
Then, the Hilbert series of $A$ is 
\begin{equation*}
    h_A(t)=5+5t+2t^2+\sum_{n\geq 3}t^n=\frac{5-3t^2-t^3}{1-t}.
\end{equation*}
\end{example}

\begin{lemma}\label{lem:h_Acond}
Suppose that $A=\kk\cQ/\cI$ is (locally) gentle. Then, 
$h_A(t^{-1})=\pm t^{k}h_A(t)$ for some $k\in\ZZ$ if and only if one of the following holds. 	
\begin{enumerate}
    \item The quiver $\cQ$ is 2-regular and $h_A(t)=\frac{\abs{\cQ_0}(1+t)}{1-t}$.
    \item The algebra $A$ is $\kk\widetilde{A}_n$ and $h_A(t)=\frac{\abs{\cQ_0}}{1-t}$.
    \item  The quiver $\cQ$ is a cycle graph with some orientation and $\cI$ is generated by all paths of length 2. Then, $h_A(t)=\abs{\cQ_0}(1+t)$.
\end{enumerate}
\end{lemma}

\begin{proof}
The fact that $h_A(t^{-1})=\pm t^{k}h_A(t)$ is clear in these three cases from the equalities $2\abs{\cQ_0}=\abs{\cQ_1}$ in (1) and  $\abs{\cQ_0}=\abs{\cQ_1}$ in (2) and (3), and that $A$ has no finite maximal paths in the first two cases and $\abs{\cQ_0}$ finite maximal paths of length 1 in the third.

Assume $h_A(t^{-1})=\pm t^{k}h_A(t)$ for some $k\in\mathbb{Z}$ and that $A$ has at least one finite maximal path. Then, we must have $k_m=\abs{\cQ_0}$ where $m$ is the length of the longest finite maximal path. Since maximal paths are disjoint, $\cQ$ must have at least $m\cdot k_m=m\abs{\cQ_0}$ arrows.

As $A$ has a finite maximal path, this implies $m<2$ (if $\abs{Q_1}=2\abs{Q_0}$, then $A$ has no finite maximal paths). Thus $m=1$, so 
\begin{equation*}
    h_A(t)=\frac{|\cQ_0|+(|\cQ_1|-|\cQ_0|)t-|\cQ_0|t^2}{1-t}.
\end{equation*}
The condition $h_A(t^{-1})=\pm t^{k}h_A(t)$ implies that $|\cQ_1|-|\cQ_0|=0$. Since $k_1=|\cQ_0|=|\cQ_1|$, every arrow must be contained in a finite maximal path of length 1. In particular, every path of length two is contained in the ideal $\cI$. 
If there exists a vertex $u$ with $\indeg(u)=2$ and $\outdeg(u)\geq 1$, then there must exist arrows $\alpha,\beta$ so that $t(\alpha)=s(\beta)=u$ and $\alpha\beta\not\in\cI$. Therefore, any vertex $u$ with $\indeg(u)=2$ is a sink, and similarly, any vertex $u'$ with $\outdeg(u')=2$ is a source.
It follows that $\indeg(v)+\outdeg(v)\leq 2$ for each $v\in\cQ_0$. Since 
\begin{equation*}
    \sum_{v\in\cQ_0}\left(\indeg(v)+\outdeg(v)\right)=2\abs{\cQ_1}
\end{equation*}
for any quiver and $|\cQ_0|=|\cQ_1|$,
we have $\indeg(v)+\outdeg(v)= 2$ for each vertex $v$. Thus, $\cQ$ must be a cycle graph with some orientation.

If $A$ has no finite maximal paths, then $h_A(t)=\frac{\abs{\cQ_0}+(\abs{\cQ_1}-\abs{\cQ_0})t}{1-t}$, so  $h_A(t^{-1})=\pm t^{k}h_A(t)$ implies either $\abs{\cQ_1}=\abs{\cQ_0}$ or  $\abs{\cQ_1}=2\abs{\cQ_0}$. The former implies that $A=\kk\widetilde{A}_n$ (see Lemma \ref{lem:invoutdeg}).
\end{proof}
  
\subsection{The center}\label{sec:center}
The center of a locally gentle algebra can be described in terms of its maximal paths. In particular, the infinite maximal paths give us generators of the following form.
\begin{defn}
If $\gamma$ is an infinite maximal path, we define the element $m_\gamma$ of $A$ as
\begin{equation*}       
    m_\gamma=\sum_{\substack{\alpha\in\cQ_1,\\\gamma_\alpha=\gamma}}\gamma_r(\alpha).
\end{equation*}
\end{defn}

Note that $m_\gamma^k=\sum_{\substack{\alpha\in\cQ_1,\\\gamma_\alpha=\gamma}}(\gamma_r(\alpha))^k$ for $k>0$ since $\gamma_r(\alpha)$ does not have any repeated arrows.

\begin{lemma}\label{lem:center} The center of $A$, $Z(A)$, is generated as a $\kk$-algebra by the following elements of $A$.
\begin{enumerate}
    \item The identity of $A$, $\one=\sum_{v\in \cQ_0}\epsilon_v$. \label{centergen1}
	\item  Finite maximal paths $\gamma \in\cB$ so that $\gamma$ is a cycle in $\cQ$. \label{centergen2}
	\item Elements $m_\gamma$ for each infinite maximal path $\gamma$ of $A$. \label{centergen3}
\end{enumerate} 
Moreover, elements of the form \eqref{centergen1}-\eqref{centergen2} and positive powers of elements of the form \eqref{centergen3} give a basis of $Z(A)$.
\end{lemma}

\begin{proof}
A quick check shows that if $p$ is an element satisfying one of \eqref{centergen1}-\eqref{centergen3}, then $\epsilon_vp=p\epsilon_v$ and $\alpha p=p\alpha$ for all $v\in\cQ_0$ and $\alpha\in\cQ_1$. It follows that $p\in Z(A)$.

Let $x=\sum_{p\in\cB} \lambda_pp$ be an element of the center. 
For each $v\in\cQ_0$, we have
\begin{equation*}
    \sum_{\substack{p\in\cB,\\s(p)=v}}\lambda_p p=\epsilon_vx=x\epsilon_v=\sum_{\substack{p\in\cB,\\t(p)=v}}\lambda_p p.
\end{equation*}
Since $\cB$ is a basis for $A$, this implies that $\lambda_p\neq 0$ if and only if $s(p)=t(p)$. Thus  $x$ is a linear combination of cycles. 

Let $\alpha$ be any arrow. The coefficient of $\alpha$ in $x\alpha$ is $\lambda_{\epsilon_{s(\alpha)}}$ while its coefficient in $\alpha x$ is $\lambda_{\epsilon_{t(\alpha)}}$, so $\lambda_{\epsilon_{s(\alpha)}}=\lambda_{\epsilon_{t(\alpha)}}$. That is, the coefficients of stationary paths are equal if there is any arrow connecting the corresponding vertices. As $\cQ$ is connected, this implies that every stationary path occurring in $x$ has the same coefficient. Hence 
$x = \lambda_\one\one+x'$, where $x'\in A_+$, $\lambda_\one\in\kk$.

If $w\in\cB_+$ with $\lambda_w\neq 0$, and $ww',w'w\in\cI$ for all $w'\in\cB_+$, then $w$ must be a finite maximal path with $s(w)=t(w)$. Namely, $w$ is an element of the form \eqref{centergen2}.

Suppose that $w\in\cB_+$ with $\lambda_w\neq 0$
and that there exists $w'\in \cB_+$ such that $ww'\notin\cI$. The coefficient of $ww'$ in $xw'$ is  $\lambda_w\neq 0$, so there exists $w''\in\cB_+$ such that $\lambda_{w''}=\lambda_w$ and $w'w''=ww'$. In particular, $F(w')=F(w)$, so $L(w)F(w)=L(w)F(w')\notin\cI$, and hence $w^2\notin\cI$. Similarly, if there exists $w'\in\cB_+$ such that $w'w\notin\cI$, then $w^2\notin\cI$. 

If $w^2\not\in\cI$, then $w$ is a cycle, and $w=(\gamma_r(\alpha_1))^k$ for some $\alpha_1\in\cQ_1$, $k>0$. Let $(\gamma_r(\alpha_1))=\alpha_1\dotsc \alpha_n$, so that $m_{\gamma_\alpha} = \sum_{i=1}^n \gamma_r(\alpha_i)$. 
Note that $\alpha_nw = (\gamma_r(\alpha_n))^k\alpha_n$. As $\alpha_nx =x\alpha_n$, the path $(\gamma_r(\alpha_n))^k$ must also occur with coefficient $\lambda_w\neq 0$. A similar argument shows that $(\gamma_r(\alpha_{i}))^k$ occurs with coefficient $\lambda_w$ for all $i$. Therefore, $\lambda_w m_{\gamma_{\alpha}}$ is a summand of $x$. The result follows.
\end{proof}

\begin{corollary}\label{cor:cyccenter} If $A$ has a unique maximal path, $\gamma$, which is infinite, then $Z(A)=\kk[m_\gamma]$ is isomorphic to the polynomial algebra $\kk[t]$ and $A$ is free $Z(A)$-module of rank $|\cQ_1|^2$. In particular, if for each vertex $v$ with out-degree 2, we choose one arrow $\alpha_v$ with source $v$, then 
$\cB_{<|\cQ_1|}\cup \setst{\gamma_r(\alpha_v)}{\outdeg(v)=2}$ is a $Z(A)$-basis for $A$.
\end{corollary}

\section{Prime ideals}\label{sec:primes}

Throughout this section, ideals are assumed to be two-sided and not necessarily graded.
\begin{defn} 
For $v\in\cQ_0$, let $\fm_v$ be the ideal generated by 
\begin{equation*}
    \setst{\epsilon_{v'}}{v'\in\cQ_0\setminus\{v\}}\cup\setst{\alpha\in\cQ_1}{s(\alpha)=t(\alpha)=v}.
\end{equation*}

Let $\gamma$ be a maximal path of $A$. The \textit{ideal corresponding to $\gamma$}, denoted $\cJ_\gamma$, is the ideal generated by all arrows with maximal path $\gamma$. That is,
\begin{equation*}
    \cJ_\gamma= \langle \left.\alpha\in\cQ_1 \right| \gamma_\alpha=\gamma\rangle =\kk\setst{w}{w\leq\gamma,\ell(w)>0} . 
\end{equation*}

Denote the left annihilator of the ideal $\cJ_\gamma$ by $\anni_A(\cJ_\gamma)$, namely $\anni_A(\cJ_\gamma)=\{x\in A: x\cJ_\gamma=0\}$.
\end{defn}

\begin{proposition}\label{prop:A'primes} Let $A=\kk\cQ/\cI$ be a locally gentle algebra with a unique maximal path, $\gamma$, which is infinite. The prime ideals of $A$ are
\begin{enumerate}
    \item the zero ideal, $0$,
    \item the ideal $\fm_v$ for each $v\in\cQ_0$, 
    \item the ideal $\langle p(m_\gamma)\rangle$ where $p(t)\in \kk[t]$ is an irreducible polynomial with $\deg(p)\geq 1$ and $p(0)\neq 0$.
\end{enumerate}
\end{proposition}

\begin{proof}
	Let $\alpha_1\dotsc\alpha_n$ be a primitive cycle of $\gamma$ and let 
    $w_{i,j}\in\cB$ be the unique path of minimal length such that $F(w_{i,j})=\alpha_i$ and $L(w_{i,j})=\alpha_j$ (so $w_{i,i}=\alpha_i$ and $w_{i,i-1}=\gamma_r(\alpha_i)=\gamma_\ell(\alpha_{i-1})$), where $\alpha_{i-kn}\defeq \alpha_i$ for any $k\in\ZZ$. 
 First, we show that ideals of the above forms are prime. 
Recall that an ideal $I\subsetneq A$ is prime if and only if $xAy\subseteq I$ implies that $x\in I$ or $y\in I$ for any $x,y\in A$ \cite{Goodearl}. 
	
\begin{enumerate}
    \item Let $x,y\in A\setminus\{0\}$. We wish to show that $xAy\not\subseteq 0$. We may choose $w_x,w_y\in\cB$ which are maximal with respect to length such that $w_x$ and $w_y$ occur with nonzero coefficient in $x$ and $y$, respectively. As $w_x$ and $w_y$ are both subpaths of $\gamma$, which is infinite and periodic, there exists $w\in\cB$ such that $w_xww_y\leq\gamma$. If $w_x'ww_y'=w_xww_y$ for some $w_x',w_y'$ occurring with nonzero coefficient in $x$ and $y$, respectively, then $w_x'$ must be an initial subpath of $w_x$ by our assumption that $w_x$ has maximal length. 
    As $w_y'$ must likewise be a terminal subpath of $w_y$, we have $w_x=w_x'$ and $w_y=w_y'$. Thus the coefficient of $w_xww_y$ in $xwy$ is the product of the coefficients of $w_x$ and $w_y$ and hence nonzero, so $0\neq xwy\in xAy$.
	
	\item For any $v\in \cQ_0$, we have $A/\fm_v\cong \kk$, so $\fm_v$ is a maximal ideal and hence prime.
	
    \item Suppose that $p\in \kk[t]$ is irreducible, $\deg(p)\geq 1$ and $p(0)=\lambda\in\kk^\times$. Let $x,y\in A$ be such that $xAy\subseteq \langle p(m_\gamma)\rangle$. Without loss of generality, we may assume that $x,y\in A_+$ (otherwise, we may replace $x$ with $x-\lambda^{-1}p(m_\gamma)x_0\in A_+$, where $x_0$ is the degree 0 portion of $x$, and likewise $y$). 
    Then there exist $p_{ij}, q_{ij}\in \kk[t]$ so that
    \begin{equation*}
        x=\sum_{1\leq i,j\leq n} p_{ij}(m_\gamma)w_{i,j} \hspace{1cm} \text{and} \hspace{1cm}y=\sum_{1\leq i,j\leq n} q_{ij}(m_\gamma)w_{i,j}.
    \end{equation*}
    Assume for contradiction that neither $x$ nor $y$ belongs to $\langle p(m_\gamma)\rangle$. Then there exist $i,j,k,l\in\{1,\ldots, n\}$ so that $p_{ij}, q_{kl}\notin \langle p\rangle\subset \kk[t]$. We will show that $ w_{i-1,l+1}\in  \langle p\rangle$ to derive a contradiction. Note that 
    $\alpha_{i-1}x\alpha_{j+1}=\alpha_{i-1}w_{i,j}  p_{ij}(m_\gamma)\alpha_{j+1}$ and likewise for $y$, so
    \begin{align*}
        \alpha_{i-1} x w_{j+1,k-1}y \alpha_{l+1}&= \alpha_{i-1} p_{ij}(m_\gamma)w_{i,j} w_{j+1,k-1} q_{kl}(m_\gamma)w_{k,l}\alpha_{l+1}\\
        &= w_{i-1,l+1}m_\gamma^s[p_{ij}q_{kl}](m_\gamma)
    \end{align*}
    for some $s\geq 0$. Thus, as $x w_{j+1,k-1}y\in\langle p(m_\gamma)\rangle$, we have  $w_{i-1,l+1}m_\gamma^s[p_{ij}q_{kl}](m_\gamma)\in \langle p(m_\gamma)\rangle$ as well.
    As $p\in\kk[t]$ is irreducible and $p(0)\neq 0$, we know that $t^sp_{ij}q_{kl}\notin \langle p\rangle $ and hence there exist polynomials $r_1, r_2\in \kk[t]$ so that $r_1p+r_2t^sp_{ij}q_{kl}=1$. Thus,
    \begin{align*}
        w_{i-1,l+1}&= w_{i-1,l+1}[r_1p+r_2t^sp_{ij}q_{kl}](m_\gamma)
        \in \langle p(m_\gamma)\rangle.
    \end{align*}
    Then $w_{i-1,l+1}=zp(m_\gamma)$ for some $z\in A$. As $p(0)\neq 0$, considering the degree zero terms yields $z\in A_+$. As $\{w_{i,j}\}_{1\leq i,j\leq n}$ is a basis for $A_+$ over $\kk[m_\gamma]$, this implies that $p(m_\gamma)z_{(i-1)(j+1)}(m_\gamma)=1$ (where $z_{(i-1)(j+1)}(m_\gamma)$ is the coefficient of $w_{i-1,j+1}$ in $z$). This contradicts our assumption that $\deg(p)\geq 1$.
\end{enumerate}

Next, we show that any nonzero proper prime ideal $\mathfrak{p}$ of $A$ must be of the form (2) or (3). We consider the cases where $\mathfrak{p}$ does and does not contain a path in $\cB$ separately.

For our first case, suppose that there exists $w\in \mathfrak{p}\cap \cB$. We will show that $\mathfrak{p}$ is of form (2).
We claim that $\mathfrak{p}$ must contain an arrow. If $\ell(w)=0$, then $w=\epsilon_v$ for some $v\in\cQ_0$, and so $\alpha=\alpha w\in \mathfrak{p}$ for some arrow $\alpha$. If $\ell(w)>1$, then $w=F(w)w'$ for some $w'\in \cB$ with $1\leq \ell(w')=\ell(w)-1$ and 
\begin{equation*}
    F(w)Aw'=\kk\{F(w)m_\gamma^kw'\}_{k\geq 0}=\kk\{w m_\gamma^k\}_{k\geq 0}\subset \langle w\rangle\subseteq \mathfrak{p}
\end{equation*}
and so either $F(w)\in \mathfrak{p}$ or $w'\in \mathfrak{p}$. As $\ell(w')<\ell(w)$, induction on $\ell(w)$ shows that $\mathfrak{p}$ contains an arrow. Without loss of generality, we may assume $\alpha_1\in\mathfrak{p}$. For any $i\neq 1$, we see that 
\begin{equation*}
    \alpha_i A\alpha_i\subseteq A\alpha_iw_{i+1,i}\subseteq \langle \alpha_1\rangle \subseteq \mathfrak{p},
\end{equation*}
since $\alpha_iw_{i+1,i}$ is the unique shortest path in $\cB$ which contains $\alpha_i$ twice and $\alpha_1\leq \alpha_iw_{i+1,i}$. Thus $\alpha_i\in \mathfrak{p}$ as well, so $\mathfrak{p}$ contains every arrow and hence $A_+\subseteq \mathfrak{p}$. Since $\mathfrak{p}\neq A$, there exists some $v_0\in\cQ_0$ so that $\epsilon_{v_0}\notin \mathfrak{p}$. For any vertex $v\neq v_0$, 
\begin{equation*}
    \epsilon_{v_0}A\epsilon_{v}\subseteq A_+\subseteq \mathfrak{p},
\end{equation*}
and hence $\epsilon_{v}\in \mathfrak{p}$. Thus $\fm_{v_0}\subseteq \mathfrak{p}$ and hence $\mathfrak{p}=\fm_{v_0}$, by maximality of $\fm_{v_0}$.

For our second case, suppose that $\mathfrak{p}$ has empty intersection with $\cB$. We will show that $wq(m_\gamma)\in \mathfrak{p}$ implies $q(m_\gamma)\in \mathfrak{p}\cap \kk[m_\gamma]$ for any polynomial $q(t)\in \kk[t]$ and $w\in\cB$. Because $m_\gamma$ is central, $wq(m_\gamma)\in \mathfrak{p}$ implies
\[wAq(m_\gamma)=w q(m_\gamma) A\subseteq \mathfrak{p}.\]
Hence, since $w\in \cB$ cannot belong to $\mathfrak{p}$, it follows that $q(m_\gamma)\in\mathfrak{p}$. We will now show that $\mathfrak{p}\cap \kk[m_\gamma]\neq 0$.  Since $\mathfrak{p}$ is nonzero, it contains a nonzero element $x$. Then, there exists some path $w$ with nonzero coefficient in $x$. Since $\gamma$ is infinite, there exist $i,j$ such that $\alpha_iw\alpha_j\notin\cI$ and hence $\alpha_ix\alpha_j\in\mathfrak{p}\setminus\{0\}$. By choice of $w_{ij}$, it follows that $\alpha_ix\alpha_j=w_{ij}q(m_\gamma)$ for some polynomial $q(t)\in \kk[t]\setminus \{0\}$. Thus, $w_{ij}\in\cB$  implies $q(m_\gamma)\in\mathfrak{p}\cap \kk[m_\gamma]$.
Thus $\mathfrak{p}\cap \kk[m_\gamma]$ is a nonzero prime ideal of $\kk[m_\gamma]$, so it is generated by $p(m_\gamma)$ for some 
irreducible monic polynomial $p(t)\in\kk[t]$. Moreover, $p(t)\neq t$, since 
$\alpha_1m_\gamma=\alpha_1\gamma_\ell(\alpha_1)\in \cB$ does not belong to $\mathfrak{p}$. Let $x\in \mathfrak{p}$. We will show that $x\in \langle p(m_\gamma)\rangle$.
Fix $v\in\cQ_0$ and take $i$ so that $s(\alpha_{i})=v$. Then 
$\alpha_{i-1}x\alpha_{i} =\alpha_{i-1}\alpha_{i}p_v(m_\gamma)$
 for some $p_v\in\kk[t]$. Since $\alpha_{i-1}x\alpha_{i}\in \mathfrak{p}$ and $\alpha_{i-1}\alpha_{i}\in\cB$, it follows that $p_v(m_\gamma)\in \mathfrak{p}\cap \kk[m_\gamma]\subseteq \langle  p(m_\gamma)\rangle$. Thus 
\begin{equation*}
    x'\defeq \sum_{v\in\cQ_0} \epsilon_vp_v(m_\gamma)\in \langle p(m_\gamma)\rangle.
\end{equation*} 
The constant coefficient of $p_v$ is exactly the coefficient of $\epsilon_v$ in $x$, so $x-x'\in A_+$ and thus
\begin{equation*}
    x-x'=\sum_{1\leq i,j\leq n}w_{i,j}p_{ij}(m_\gamma)
\end{equation*}
for some $p_{ij}\in \kk[t]$. As $x-x'\in \mathfrak{p}$, we have
\begin{equation*}
    \alpha_{i-1}(x-x')\alpha_{j+1}=\alpha_{i-1}w_{i,j}\alpha_{j+1}p_{ij}(m_\gamma)\in \mathfrak{p}
\end{equation*}
for any $1\leq i,j\leq n$. Noting that $\alpha_{i-1}w_{i,j}\alpha_{j+1}\in\cB$, we have $p_{ij}(m_\gamma)\in \mathfrak{p}\cap\kk[m_\gamma]\subseteq \langle  p(m_\gamma)\rangle$. Thus $x=x'+\sum_{1\leq i,j\leq n}w_{i,j}p_{ij}(m_\gamma)\in  \langle p(m_\gamma)\rangle$ as well. Hence, $\mathfrak{p}= \langle p(m_\gamma)\rangle$.
\end{proof}

\begin{corollary}\label{cor:annihilator} Suppose that $A$ has an infinite maximal path $\gamma$. Then, $\anni_A(\cJ_\gamma)$ is a  two-sided ideal of $A$ and  $A=A'\oplus \anni_A(\cJ_\gamma)$  where $A'$ is a locally gentle subalgebra of $A$ with unique maximal path $\gamma$.
\end{corollary}
\begin{proof} 		
Let $\cQ'$ be the minimal subgraph of $\cQ$ containing $\gamma$. That is, $\cQ'_1=\setst{\alpha\in\cQ_1}{\gamma_\alpha=\gamma}$ and $\cQ_0'=\setst{s(\alpha)}{\alpha\in \cQ'_1}$ (since $\gamma$ is infinite, this includes every vertex $\gamma$ passes through). Then $A'\defeq  \kk\cQ'/(\cI\cap \kk\cQ')$ is locally gentle and has a unique maximal path, $\gamma$. View $A'$ as a subalgebra of $A$, so that $A'_+= \cJ_\gamma$. By Proposition \ref{prop:A'primes}, $0$ is a prime ideal of $A'$ and so $xA'\alpha\neq 0$ for any $x\in A'\setminus \{0\}$ and $\alpha\in \cQ'_1$. Since $\sum_{\alpha\in\cQ_1'}A'\alpha =  A'_+=\cJ_\gamma$, this shows that $\anni_A(\cJ_\gamma)\cap A'=0$.

As $\gamma$ is maximal and $\cJ_\gamma$ has a basis of non-stationary subpaths of $\gamma$, we see that $\anni_A(\cJ_\gamma)$ contains any path which is not a subpath of $\gamma$. That is, the paths in $\overline{\cB}\defeq \cB\setminus A'$. Thus $A=A'\oplus \anni_A(\cJ_\gamma)$. Moreover, since $\anni_A(\cJ_\gamma)\supseteq \spn \overline{\cB}$ and $A'\oplus  \spn \overline{\cB}=A$, it must be that $\anni_A(\cJ_\gamma)= \spn \overline{\cB}$, which 
 is a two-sided ideal of $A$.
\end{proof}
\begin{theorem}\label{thm:primeideals} Let $A=\kk\cQ/\cI$ be (locally) gentle. 	
The prime ideals of $A$ are exactly
\begin{enumerate}
    \item the ideal $\fm_v$ for each $v\in\cQ_0$,
    \item the ideal $\anni_A(\cJ_\gamma)$ for any infinite maximal path $\gamma$ in $A$.
    \item  the ideal 
    \begin{equation*}
        \mathfrak{p}(\gamma,p)\defeq \anni_A(\cJ_\gamma)+\langle p(m_\gamma)\rangle, 
    \end{equation*} 
    for an infinite maximal path $\gamma$ and an irreducible polynomial $p(t)\in \kk[t]$ such that  $\deg(p)\geq 1$ and $p(0)\neq 0$.
\end{enumerate}
	
In particular, $A$ is prime (i.e., $0$ is a prime ideal) if and only if $A$ has only one maximal path and that path has infinite length.  The height zero prime ideals of $A$ are the ideals  $\anni_A(\cJ_\gamma)$ for $\gamma$ maximal and $\fm_v$ for those vertices $v$ which no infinite maximal path passes through. If $\gamma$ is an infinite maximal path in $A$ which passes through $v$ and $p\neq t$ is an irreducible monic polynomial, the inclusions are as follows.
\begin{center}
    \begin{tikzpicture}
        \node at (0,-.25) {$\anni_A(\cJ_\gamma)$};
        \node at (-1,1.25) {$\fm_v$};
        \node at (1,1.25) {$\mathfrak{p}(\gamma,p)$};
        \draw[-] (-1,1) to (-.25,.2);
        \draw[-] (1,1) to (.25,.2);
    \end{tikzpicture}
\end{center}
\end{theorem}
 
\begin{proof}	
	
Let $\mathfrak{p}$ be a proper prime ideal of $A$. Suppose $\mathfrak{p}$ contains every arrow and hence $A_+\subseteq \mathfrak{p}$. As $\mathfrak{p}$ is proper, there exists some vertex $v\in \cQ_0$ such that $\epsilon_v\notin \mathfrak{p}$. For any other vertex $v\neq v'\in \cQ_0$, we have $\epsilon_vA\epsilon_{v'}\subseteq A_+$, so we must have $\epsilon_{v'}\in \mathfrak{p}$ as well. Thus, 
\begin{equation*}
    \fm_v=\langle \epsilon_{v'}, a\rangle_{\substack{v'\in\cQ_0\setminus \{v\},\\ a\in \cQ_1}}\subseteq \mathfrak{p}.
\end{equation*}
As $\fm_v$ is maximal, $\mathfrak{p}=\fm_v$.
Otherwise, there exists an arrow $a_0\in\cQ_1$ such that $a_0\notin \mathfrak{p}$. Then $a_0Aa_0\not\subseteq\mathfrak{p}$, so $a_0Aa_0\neq 0$. Thus $\gamma\defeq \gamma_{a_0}$ is infinite.
As
$Aa_0\subset \cJ_\gamma$, we see that $xAa_0=0\subseteq \mathfrak{p}$ for any $x\in\anni_A(\cJ_\gamma)$ and hence $\anni_A(\cJ_\gamma)\subseteq \mathfrak{p}$.	
As noted in Corollary \ref{cor:annihilator}, $A=A'\oplus \anni_A(\cJ_\gamma)$, where $A'$ is the  subalgebra of $A$ with unique maximal path $\gamma$. Thus $\mathfrak{p}= \mathfrak{p}\cap A'+\anni_A(\cJ_\gamma)$ is prime if and only if $\mathfrak{p}\cap A'$ is prime in $A'$, so the result follows from Proposition \ref{prop:A'primes}.
\end{proof}

\begin{example}
Let $\kk$ be algebraically closed. 
The prime ideals of $\kk[x,y]/\langle xy,yx\rangle$ are $\fm_1=\langle x,y\rangle$,  $\anni(\cJ_{\gamma_x})=\langle y\rangle$, $\anni(\cJ_{\gamma_y})=\langle x\rangle$, $\mathfrak{p}(\gamma_x,t-\lambda)=\langle x-\lambda, y\rangle $ and $\mathfrak{p}(\gamma_y, t-\lambda)=\langle y-\lambda, x\rangle $ for $\lambda\in \kk^\times$. 
\begin{center}
    \begin{tikzpicture}
        \node at (0,-.25) {$\anni(\cJ_{\gamma_x})$};
        \node at (1.25,1.25) {$\fm_1$};
        \node at (-1,1.25) {$\mathfrak{p}(\gamma_x, t-\lambda)$};
        \draw[-] (-1,1) to (-.25,.2);
        \draw[-] (1,1) to (.25,.2);
        \node at (2.5,-.25) {$\anni(\cJ_{\gamma_y})$};
        \node at (3.5,1.25) {$\mathfrak{p}(\gamma_y, t-\lambda)$};
        \draw[-] (1.5,1) to (2.25,.2);
        \draw[-] (3.5,1) to (2.75,.2);
    \end{tikzpicture}
\end{center}
 
Let $A=\kk\cQ/\langle ab,cd\rangle$ where $\cQ$ is the quiver below.
\begin{center} 
    \begin{tikzpicture}[scale=1]
        \vtx{0,0}{1}
        \vtx{1,0}{2}
        \vtx{2,0}{3}
        \draw[->] (.15,0) to node[above]{$b$}(.85,0) ;
        \draw[->] (1.15,0) to node[above]{$c$}(1.85,0) ;
        \draw[->] (-.15,.1) to[out=150, in=-150, looseness=10] node[left]{$a$} (-.15,-.1) ;
        \draw[->] (2.15,.1) to[out=30, in=-30, looseness=10] node[right]{$d$} (2.15,-.1) ;
    \end{tikzpicture}
\end{center}
The prime ideals of $A$ are $\fm_1=\langle \epsilon_2,\epsilon_3, a\rangle $, $\fm_2=\langle \epsilon_2,\epsilon_3\rangle $, $\fm_3=\langle \epsilon_1,\epsilon_2, d\rangle $, $\anni_A(\cJ_{\gamma_a})=\langle \epsilon_2, \epsilon_3\rangle$, 
 $\anni_A(\cJ_{\gamma_d})=\langle \epsilon_1, \epsilon_2\rangle$, $\mathfrak{p}(\gamma_a, t-\lambda)=\langle a-\lambda, \epsilon_2, \epsilon_3\rangle $, and 
$\mathfrak{p}(\gamma_d, t-\lambda)=\langle d-\lambda, \epsilon_1, \epsilon_2\rangle $ for $\lambda\in\kk^\times$. 
\begin{center}
    \begin{tikzpicture}
        \node at (-.5,-.25) {$\anni_A(\cJ_{\gamma_a})$};
        \node at (-1.5,1.25) {$\mathfrak{p}(\gamma_a,t-\lambda)$};
        \node at (.5,1.25) {$\fm_1$};
        \draw[-] (-1.5,1) to (-.75,.2);
        \draw[-] (.5,1) to (-.25,.2);
        \node at (2,1.25) {$\fm_2$};
        \node at (4.5,-.25) {$\anni_A(\cJ_{\gamma_d})$};
        \node at (3.5,1.25) {$\fm_3$};
        \node at (5.5,1.25) {$\mathfrak{p}(\gamma_d,t-\lambda)$};
        \draw[-] (3.5,1) to (4.25,.2);
        \draw[-] (5.5,1) to (4.75,.2);
    \end{tikzpicture}
\end{center}
\end{example}

\begin{corollary}\label{cor:semiprime}
The prime radical of $A$ is $\langle \alpha\in\cQ_1\mid\gamma_\alpha \text{ is finite}\rangle.$
In particular, $A$ is semiprime if and only if it has  no finite maximal paths.
\end{corollary}

\begin{proof}
The prime radical of $A$ is 
\begin{equation*}
    \bigcap_{\mathfrak{p}\in\text{Spec}(A)}= \bigcap_{v\in\cQ_0}\fm_v\cap \bigcap_{\gamma\text{ infinite maximal}} \anni_A(\cJ_\gamma).
\end{equation*}
It is clear that $\bigcap_{v\in\cQ_0}\fm_v=A_+$ while
\begin{align*}
    \bigcap_{\gamma\text{ infinite maximal}} \anni_A(\cJ_\gamma)&= \bigcap_{\gamma\text{ infinite maximal}} \spn\setst{w\in\cB}{ w\text{ is not a subpath of }\gamma}\\& = \spn\setst{w\in\cB}{ w\text{ is not a subpath of any infinite maximal path}}
\end{align*}
by Corollary \ref{cor:annihilator}. Thus the prime radical is spanned by all non-stationary subpaths of finite maximal paths. Such paths are
generated by the arrows contained in finite maximal paths, so
\begin{equation*}
    \bigcap_{\mathfrak{p}\in\text{Spec}(A)} \mathfrak{p}=\langle \alpha\in\cQ_1\mid\gamma_\alpha \text{ is finite}\rangle.
\end{equation*}
In particular, this intersection is $A_+$ if $A$ is gentle and $0$ if $A$ has no finite maximal paths. As $A$ is semiprime when the prime radical is 0 \cite{Goodearl}, the result follows.
\end{proof}
 
\section{Homological dimensions}\label{sec:homological-dim}
 
In \cite[Theorems 5.10 and 6.9]{LGH}, it was shown that for a gentle algebra, the global and injective dimension are determined by the maximal paths of Koszul dual. In this section, we will extend this description to locally gentle algebras. In Section \ref{subsec:global-dim}, we construct minimal graded projective resolutions for the graded simple modules and use this to find the graded global dimension of $A$. In Section \ref{subsec:injdim}, we construct minimal graded injective resolutions of the indecomposable graded projective modules and use this to find the graded injective dimension of $A$. The reader may consult \cite{toolkit} for a nice overview of the definitions of these concepts in the graded case. 

We construct our resolutions in the category $\grmod{A}$ of graded right $A$-modules, where $A=\kk\cQ/\cI$ is a (locally) gentle algebra as in Notation \ref{not:algebra}.
For the remainder of this section, we fix a vertex $v\in\cQ_0$.

\subsection{Global dimension} \label{subsec:global-dim}

\begin{notation}
We denote the graded simple module of $A$ at the vertex $u \in \cQ_0$ by $S(u)$. Note that this is the module $\kk\set{\epsilon_u}$ where the right action of $A$ is given by $\epsilon_u\cdot \epsilon_u=\epsilon_u$ and $\epsilon_u\cdot x=0$ for any path $x\in \cB$ so that $x\neq \epsilon_u$.
\end{notation}

\begin{defn}
Our projective resolution traces the forbidden right maximal paths beginning with arrows in $s^{-1}(v)$. Thus, it is helpful to define the set of forbidden paths
\begin{equation*}
    \Rv\defeq \setst{p\in \cB^\#}{s(p)=v} =\setst{p}{p\leq(\gamma^\#_r(\alpha))^k,\alpha\in s^{-1}(v),s(p)=v,k>0}.
\end{equation*} 
Note that these forbidden paths are exactly the initial subpaths of $(\gamma_r^\#(\alpha))^k$ where $\alpha\in s^{-1}(v)$ and $k>0$.
Let $\cR'(v)\subseteq \Rv$ denote the set of right maximal forbidden paths with source $v$.
For a path $p \in \Rv_k$ and $k>0$, define 
\begin{equation*}
    p_*: (\epsilon_{t(p)}A)[-k]\to (\epsilon_{s(L(p))}A)[-k+1], \quad\text{where} \quad p_*(x)=L(p)x,
\end{equation*}
and $L(p)$ is an arrow in $A$. We define the resolution $P_\bullet\to S(v)$, given by
\begin{equation}\label{eq:proj-res}
\begin{tikzcd}
   \dotsc\arrow{r}{d_{k+1}}&P_k\arrow{r}{d_k}&P_{k-1}\arrow{r}{d_{k-1}}&\dotsc\arrow{r}{d_1}&P_0\arrow{r}{\pi}&S(v)\arrow{r}&0,
\end{tikzcd}
\end{equation}
where for $i>0,k>1$, we let
\begin{equation*}
    P_0=\epsilon_v A, \quad P_i=\bigoplus_{p\in \Rv_i} \epsilon_{t(p)} A[-i], \quad d_1=\sum_{p\in\Rv_1} p_*, \quad
    d_{k}=\bigoplus_{p\in \Rv_k}p_*,
\end{equation*}
and $\pi$ is the natural projection.
\end{defn}

\begin{rmk}\label{rmk:Rvinduct}
Note that $\abs{\Rv_k}\leq 2$ for each $k>0$, and $\Rv_0=\set{\epsilon_v}$. If $p\in \Rv_k$, then there exists at most one arrow $p_{k+1}$ in $A$ so that $L(p)p_{k+1}\in\cI$. If such an arrow $p_{k+1}$ exists, then $pp_{k+1}\in \Rv_{k+1}$. If no such arrow exists, then $\abs{\Rv_{k+1}}<\abs{\Rv_k}$ and $p$ is a right maximal forbidden path of length $k$. More precisely, we have
\begin{equation*}
    \Rv_{k+1}=
    \setst{p\alpha}{{p\in \Rv_k,\alpha\in s^{-1}(t(p)), L(p)\alpha\in\cI}}.
\end{equation*}
Note that if $p,p'\in \Rv_k$, then $L(p)=L(p')$ implies $p=p'$ (as $A^\#$ is gentle and these paths have the same length).
\end{rmk}

\begin{lemma}\label{lem:minprojres}
Let $A=\kk\cQ/\cI$ be a (locally) gentle algebra, and fix $v\in\cQ_0$. The resolution in \eqref{eq:proj-res} is a minimal graded projective resolution of $S(v)$.
\end{lemma}

\begin{proof}
It is clear that each $P_i$ is projective since $A=\bigoplus_{v\in\cQ_0}\epsilon_vA$. For $k>1$, let $p\in \Rv_k$. Then, we have that $\im(p_*)=L(p)A$. There is at most one arrow $p_{k+1}\in\cQ_1$ so that $s(p_{k+1})=t(p)$ and $L(p)p_{k+1}\in\cI$. Then, if $p_{k+1}$ exists, it is clear that $p_{k+1}A \subseteq \Ker(p_*)$.
If $q\in\Ker(p_*)$ is a path of positive length, then $L(p) q \in \cI$, and since $q\not\in\cI$, it must be the case that $p_{k+1}=F(q)$. Therefore, we have that $\Ker{(p_*)}=p_{k+1}A$. 
It follows that $\Ker{(p_*)}=p_{k+1} A=\im((pp_{k+1})_*)$. 
If no arrow $p_{k+1}$ exists, then $pq\not\in\cI$ for any $q\in\epsilon_{t(p)}A$, so $\Ker(p_*)$ is trivial. 
By Remark \ref{rmk:Rvinduct}, we have that 
\begin{equation*}
    \im(d_k)=\bigoplus_{p\in \Rv_k}L(p)A \quad\text{and}\quad \Ker(d_k)=\bigoplus_{\substack{p\in \Rv_k,\\ \alpha\in s^{-1}(t(p)),\\L(p)\alpha\in\cI}}\alpha A=\bigoplus_{p\in \Rv_{k+1}}L(p)A.
\end{equation*}
Therefore, $\im(d_{k+1})=\Ker(d_k)$ for $k>1$.

We have that $\im(d_1)=\sum_{\alpha\in s^{-1}(v)}\alpha A$, which is precisely submodule of $P_0$ generated by all paths of positive length i.e. $\Ker(\pi)$.
Further, if $v$ is the source of two distinct arrows $\alpha_1$ and $\beta_1$, then $d_1((x,y))=\alpha_1x+\beta_1y\in\cI$ for $x\in \epsilon_{t(\alpha_1)}A, y\in \epsilon_{t(\beta_1)}A$ it must be the case that both $\alpha_1x\in\cI$ and $\beta_1y\in\cI$. If $\alpha_1x\in\cI$, then $x$ must either be 0 or have positive length with $F(x)=\alpha_2$ for $\alpha_1\alpha_2\leq\gamma_r^\#(\alpha_1)$. A similar argument for $y$ shows that $\Ker(d_1)\subseteq \im(d_2)$.

To show this resolution is minimal, we show that  each map $d_n$ and $\pi$ are graded projective covers. In other words, we show that $\Ker(d_n)\subseteq P_n\grJ(A)$ for all $n\geq 1$ and $\Ker(\pi)\subseteq P_0\grJ(A)$.  Note that for a vertex $w$, $\epsilon_wA \grJ(A)=\bigoplus_{\alpha\in s^{-1}(w)}\alpha A$. Therefore, we have that
\begin{equation*}
    P_k\grJ(A) = \bigoplus_{p\in \Rv_k} \epsilon_{t(p)} A[-k]\grJ(A)=\bigoplus_{\substack{\alpha\in s^{-1}(t(p)),\\p\in \Rv_k}}\alpha A,
\end{equation*}
and it follows immediately that $\Ker(d_k)\subseteq P_k\grJ(A)$ for $k>1$. Further, $P_0\grJ(A)=\epsilon_vA\grJ(A)=\bigoplus_{\alpha\in s^{-1}(v)}\alpha A=\Ker(\pi)$. Therefore, this is a minimal graded projective resolution.
\end{proof}

\begin{example}
Let $A$ be as in Example \ref{ex:hilb}. Then, $\cR(3)=\{\alpha_1, \alpha_1c, \alpha_1c\alpha_2, \beta_1, \beta_1\beta_2\}$, and a minimal projective resolution of $S(3)$ is
\begin{align*}
    \epsilon_1A[-3]\xrightarrow{\left(\begin{matrix}(\alpha_1c\alpha_2)_*\\0\end{matrix}\right)} \epsilon_2A[-2]\oplus \epsilon_5A[-2]\xrightarrow{\left(\begin{matrix}(\alpha_1c)_*& 0\\0&(\beta_1\beta_2)_*\end{matrix}\right)} \epsilon_2A[-1]\oplus  \epsilon_4A[-1] \\\xrightarrow{\left(\begin{matrix}(\alpha_1)_*&(\beta_1)_*\end{matrix}\right)}\epsilon_3A\rightarrow S(3).
\end{align*}
Thus, $\pdim(S(3))= \max\set{\ell(\alpha_1c\alpha_2), \ell(\beta_1\beta_2)}=3$.
\end{example}

The beginning of this minimal projective resolution matches the one described in \cite[Lemma 3.7]{RR19}. The following theorem extends the result of \cite[Theorem 5.10]{LGH} to locally gentle algebras.

\begin{theorem} \label{thm:gldim}
The global dimension of $A$  is the maximum of the lengths of its maximal forbidden paths. That is, it is the length of the longest maximal path in $A^\#$ (where the length of an infinite maximal path is considered to be $\infty$).
\end{theorem}

\begin{proof}
Let $S\defeq A/\grJ(A)=A_0$. As $S$ is separable, 
the global dimension of $A$ is exactly the projective dimension of $S$ \cite[Proposition 3.18]{toolkit}. Since $S\cong\bigoplus_{v\in\cQ_0}S(v)$ as right $A$-modules,  
\begin{equation*}
    \lgldim(A)=\rgldim(A)=\pdim(S_A)=\max_{v\in\cQ_0}\pdim(S(v)).
\end{equation*}
This is the length of the longest maximal path $\gamma$ in $A^\#$: If $\gamma$ is finite, then $\pdim_A(S(s(\gamma))=\ell(\gamma)$ and $\pdim_A(S(v))\leq \ell(\gamma)$ for any other vertex $v$. Otherwise $\gamma$ is infinite and $\pdim_A(S(v))=\infty$ for any vertex $v$ which $\gamma$ passes through. 
\end{proof}

Now, we use the projective resolution found above to make some explicit calculations of $\Ext^i(S(v),A)$.
\begin{defn}
For a vertex $u$ and $x\in A\epsilon_u$ of degree $k$, we let $\phi_x:\epsilon_uA[-k]\to A$ denote the homomorphism which sends $\epsilon_u$ to $x$, so $\phi_x(p)=xp$ for any $p\in \epsilon_uA[-k]$. Similarly, if $u'\in\cQ_0$ and $(x,y)\in \left(A\epsilon_{u}\oplus A\epsilon_{u'}\right)$ of degree $k$, we let $\phi_{(x,y)}$ denote the homomorphism in $\Hom_A(\left(\epsilon_{u}A\oplus \epsilon_{u'}A\right)[-k],A)$ given by $\phi_{(x,y)}(p,p')=xp+yp'$.
\end{defn}

There is the one-to-one correspondence between the elements of degree $i$ in $\Hom_A (\epsilon_uA[-k],A)$ and elements of degree $i+k$ in $A\epsilon_u$, given by $\phi\mapsto\phi(\epsilon_u)$.  
In particular, 
\begin{equation*}
    \Hom^i_A(\epsilon_{u}A[-k],A)= \setst{\phi_x}{x\in (A\epsilon_{u})_{i+k}}\cong (A\epsilon_u)_{i+k},
\end{equation*}
which yields the following isomorphism of graded left $A$-modules
\begin{equation*}
    \grHom_A(\epsilon_{u}A[-k],A)=\bigoplus_{i\geq -k} \Hom^i_A(\epsilon_{u}A[-k],A)\cong \bigoplus_{i\geq -k} (A\epsilon_u)_{i+k}=(A\epsilon_u)[k].
\end{equation*}
 Then, $\grHom_A(\left(\epsilon_{u}A\oplus \epsilon_{u'}A\right)[-k],A) \cong (A\epsilon_{u}\oplus A\epsilon_{u'})[k]$. It follows that 
\begin{equation*}
    \grHom_A(P_k,A)\cong\bigoplus_{p\in\Rv_k}(A\epsilon_{t(p)})[k].
\end{equation*}
Let $d_k^*:\grHom_A(P_{k-1},A)\to\grHom_A(P_k,A)$ be the map induced by $d_k$, namely $d_k^*(\phi)(y)=\phi(d_k(y))$ for $y\in P_k$ and $\phi\in\grHom_A(P_{k-1},A)$. Note that, if $p\in R(v)_{k}$ and  $p'\in R(v)_{k-1}$, then for $\phi_x\in \grHom_A(\epsilon_{t(p')}A[-(k-1)],A)$ we have 
 $(p_*)^*(\phi_x)=\phi_{xd_k(\epsilon_{t(p)})}=\phi_{xL(p)}$
 if $p=p'L(p)$ and $0$ otherwise.

\begin{proposition}\label{prop:extcalc} Keep the hypotheses of Lemma \ref{lem:minprojres}. We have the following isomorphisms of graded left $A$-modules.
\begin{enumerate}
    \item $\grExt^0_A(S(v),A)\cong 
    \begin{cases}
    A\epsilon_v & \text{if }\outdeg(v)=0,\\
    A\alpha & \text{if }\outdeg(v)=1,\alpha\in t^{-1}(v),\alpha'\in s^{-1}(v),\alpha\alpha'\in\cI,\\
    0 & \text{otherwise}.
    \end{cases}$
    
    \item $\grExt^1_A(S(v),A)\cong 
    \begin{cases}
    \left(A\epsilon_{t(\alpha)}/A\alpha\right)[1] & \text{if } s^{-1}(v)=\{\alpha\} = \cR'(v)_1,\\
    S(v) & \text{if }\outdeg(v)=2, \cR'(v)_1=\emptyset,\\
    0 & \text{if }\outdeg(v)<2, \cR'(v)_1=\emptyset.\\
    \end{cases}$
    
    The case where $\outdeg(v)=2$ and $\cR'(v)_1\neq \emptyset$ is not shown here, but $\grExt^1_A(S(v),A)$ is also nonzero in that case.
    
    \item For $i>1$,
    \begin{equation*}
        \grExt^i_A(S(v),A)\cong \bigoplus_{p\in \cR'(v)_i}\left( A\epsilon_{t(p)}/AL(p)\right)[i].
    \end{equation*}
\end{enumerate}
\end{proposition}

\begin{proof}
(1) To see the conditions on $\grExt_A^0(S(v),A)$, let $\phi\in \grHom_A(S(v),A)$.
Then $\phi(\epsilon_v)x=\phi(\epsilon_vx)=0$ for every $x\in A_+$ and hence $\phi(\epsilon_v)$ is a sum of right maximal paths with target $v$. If $w$ is not right maximal, then there is an arrow $\alpha$ so that $w\alpha\notin\cI$ and $w$ has the same coefficient in  $\phi(\epsilon_v)$ as $w\alpha$ has in $\phi(\epsilon_v)\alpha = 0$. 

If $\outdeg(v)=0$, then any path which terminates at $v$ is right maximal and so $\grHom_A(S(v),A)\cong A\epsilon_v$. If $\outdeg(v)=2$, then no right maximal path terminates at $v$  and so $\grHom_A(S(v),A)=0$.

Assume that $\outdeg(v)=1$ and that there is a right maximal path with target $v$. Let $\alpha'$ denote the arrow with source $v$. Since there is a right maximal path with target $v$, there exists an arrow $\alpha$ with target $v$ such that $\alpha\alpha'\in\cI$ and a path $p\in\cB$ which terminates at $v$ is right maximal if and only if $L(p)=\alpha$. Thus $\grHom_A(S(v),A)\cong A\alpha$.

(2) As $\cR(v)_1=s^{-1}(v)$, we have $\grExt^1_A(S(v),A)=0$ when $\outdeg(v)=0$. If $s^{-1}(v)=\set{\alpha}$, then
\begin{equation*}
    \im (d^*_1)= \setst{\phi_{x}}{x\in A\alpha}\subseteq \grHom_A\left((A\epsilon_{t(\alpha)})[-1],A\right)=\grHom_A(P_1,A).
\end{equation*}
If $s^{-1}=\set{\alpha,\beta}$ ($\outdeg(v)=2$), we instead have
\begin{equation*}
    \im (d^*_1)
    = \setst{\phi_{(x\alpha, x\beta)}}{x\in A\epsilon_v}\subsetneq \setst{\phi_x}{x\in A\alpha\oplus A\beta}=\setst{\phi_{x}}{x\in \bigoplus_{p\in \Rv_1}AL(p)}.
\end{equation*} 
In both cases,
\begin{equation*}
    \Ker (d^*_{2}) = \setst{\phi_{x}}{ x\in \bigoplus_{\alpha\in s^{-1}(v)}Ax_\alpha}, \quad \text{where}\quad
    x_\alpha=\begin{cases}
        \epsilon_{t(\alpha)} & \text{if }\alpha\in \cR'(v),\\
        \alpha & \text{otherwise.}
    \end{cases}
\end{equation*}
Clearly, $\grExt^1_A(S(v),A)\neq 0$ if $\Rv_1'\neq \emptyset$ and $\grExt^1_A(S(v),A)=0$ if $\outdeg(v)=1$ and $\Rv_1'=\emptyset$. 

Suppose that $s^{-1}(v)=\{\alpha,\beta\}$ and $\cR'(v)_1=\emptyset$. 
Let $(x,y)\in (A\alpha\oplus A\beta)_{\geq 2}$. For any $w\in\cB_+$ with $t(w)=v$, we have $w\alpha\notin\cI$ if and only if $w\beta\in\cI$. Hence there exist $x',y'\in A_+\epsilon_v$ such that $(x,y)=(x'\alpha, y'\beta)$ and $x'\beta,y'\alpha\in\cI$. Thus
\begin{equation*}
    \phi_{(x,y)}=\phi_{((x'+y')\alpha,(x'+y')\beta)}=d_1^*(\phi_{x'+y'})\in \im(d_1^*).
\end{equation*}
So we see that $\Ker(d^*_2)_{\geq 1}=\setst{\phi_{(x,y)}}{(x,y)\in   ((A\alpha \oplus A\beta)[1])_{\geq 1}}=\im(d^*_1)_{\geq 1}$, and hence
\begin{equation*}
    \grExt_A^1(S(v),A)\cong \left((A\alpha\oplus A\beta)[1]\right)_0/\kk\{(\alpha,\beta)\}\cong S(v).
\end{equation*}
as graded left $A$-modules.

(3) When $i>1$, we have
\begin{equation*}
    \im (d^*_i)= \setst{\phi_{x}}{x\in \bigoplus_{p\in \Rv_i}AL(p)}\subseteq \bigoplus_{p\in \Rv_i} \grHom_A\left((A\epsilon_{t(p)})[-i],A\right)=\grHom_A(P_i,A),
\end{equation*} 
and
\begin{equation*}
    \Ker (d^*_{i+1}) = \setst{\phi_{x}}{ x\in \bigoplus_{p\in \cR(v)_i}Ax_p}, \quad \text{where}\quad
    x_p=\begin{cases}
        \epsilon_{t(p)} & \text{if }p\in \cR'(v),\\
        L(p) & \text{otherwise.}
    \end{cases}
\end{equation*}
The result follows. 
\end{proof}

The following is an immediate consequence of the above proposition.

\begin{corollary}\label{cor:extcalc}
Let $A=\kk\cQ/\cI$ be (locally) gentle and let $v\in \cQ_0$. 
Then, $\grExt^i_A(S(v),A)\neq 0$ if and only if one of the following holds.
\begin{enumerate}
    \item There is at least one finite maximal path with target $v$ and $i=0$.
    \item There are two arrows with source $v$ and $i=1$.
    \item There is a right maximal forbidden path with source $v$ of length $i>0$.
\end{enumerate}
\end{corollary}

\begin{corollary}\label{cor:depth}
    Let $A=\kk\cQ/\cI$ be (locally) gentle. Then $A$ has depth 0 if and only if it has a finite maximal path and depth 1 otherwise.
\end{corollary}
\begin{proof}
Since $\grJ(A)=A_+$, we know that $S=A/\grJ(A)$ is isomorphic to $\bigoplus_{v\in\cQ_0}S(v)$ as $(A,A)$-bimodules. As $\cQ$ is finite,
\begin{equation*}
    \dpth(A)=\min\setst{i}{ \grExt_A^i(S(v),A)\neq 0\text{ for some }v\in\cQ_0}.
\end{equation*}
Therefore, by Corollary \ref{cor:extcalc}, we know that $\dpth(A)=0$ if and only if $A$ has a finite maximal path.  By Lemma \ref{lem:invoutdeg}, if $A$ has no finite maximal path, then either $A=\kk\tilde{A}_n$ or $\cQ$ has a vertex which is the source of two arrows, which implies that $\dpth(A)\leq 1$. As every maximal forbidden path of $\kk\tilde{A}_n$ has length one, we know that $\dpth(A)\leq 1$ as well.
\end{proof}

\subsection{Injective dimension} \label{subsec:injdim}

In this subsection, we will construct a graded injective resolution of the graded projective module $\epsilon_vA$. Similar to the minimal graded projective resolution of $S(v)$ constructed in Lemma \ref{lem:minprojres}, the graded injective resolution corresponds to a set of forbidden paths.

\begin{defn} 
Let $\cL'$ denote the set of all finite maximal forbidden paths. For $p\in \cL'_n$ and $1\leq i\leq n$, we denote the $(n+1-i)$th arrow in $p$ by $p_{-i}$, so that $p=p_{-n}\dotsc p_{-2}p_{-1}$.

Let $\Lplus$ be the set of pairs $(p,w)\in \cL'\times \cB$ such that $s(w)=v$ and $t(w)=t(p)$ but $L(w)\neq p_{-1}$. 
That is,
\begin{equation*}
    \Lplus=\setst{(p,w)\in \cL'\times \cB}{t(p)=t(w),L(p)\neq L(w), s(w)=v}.
\end{equation*} 
\end{defn}

For $(p,w)\in \Lplus$, we have $L(p)\alpha\notin \cI$ for any arrow $\alpha$ with target $t(p)=t(w)$. Consequently, as $L(w)\neq L(p)$, either $w=\epsilon_v$ or $w$ is right maximal. In the latter case, $w=\gamma_r(\alpha)$ for some $\alpha\in s^{-1}(v)$ with $\gamma_\alpha$ rfinite. The following example illustrates the two varieties of pairings.
\begin{example}
  Suppose $\cQ$ contains the subgraph
\begin{center}
    \begin{tikzpicture}
    
    \vtx{-1.5,0}{}
    \vtx{0,0}{}
    
    \vtx{1.5,0}{}
    \vtx{3,0}{}
    \vtx{5,0}{}
    \vtx{6.5,0}{$v$}
    \vtx{8,0}{}
    \draw[->] (.2,0) to node[midway, above]{\footnotesize{$p_{-n}$}}(1.3,0) ;
    \draw[->] (1.7,0) to node[midway, above]{\footnotesize{$p_{-(n-1)}$}}(2.8,0) ;
    \draw[->] (3.2,0) to (3.7,0) ;
    \node at (4,0) {$\cdots$};
    \draw[->] (4.3,0) to (4.8,0) ;
    \draw[->] (5.2,0) to node[midway, above]{\footnotesize{$p_{-1}$}}(6.3,0) ;
    \draw[->, dotted] (6.7,0) to node[midway, above]{\footnotesize{$\beta_2$}}(7.8,0) ;
    \draw[->, dotted] (-1.3,0) to node[midway, above]{\footnotesize{$\beta_{1}$}}(-.2,0) ;
\end{tikzpicture}

\begin{tikzpicture}
    \vtx{-1.5,0}{}
    \vtx{0,0}{}
    \vtx{1.5,0}{}
    \vtx{3,0}{}
    \vtx{5,0}{}
    \vtx{6.5,0}{}
    \vtx{8,0}{}
    \vtx{1.5,0}{}
     \vtx{1.5,-.8}{$v'$}
     \vtx{3,-.8}{}
    \vtx{5,-.8}{}
    \vtx{6.5,0}{}
    \draw[->,dotted] (-1.3,0) to node[midway, above]{\footnotesize{$\beta'_{1}$}}(-.2,0) ;
    \draw[->] (.2,0) to node[midway, above]{\footnotesize{$p'_{-n'}$}}(1.3,0) ;
    \draw[->] (1.7,0) to node[midway, above]{\footnotesize{$p'_{-(n'-1)}$}}(2.8,0) ;
    \draw[->] (3.2,0) to (3.7,0) ;
    \node at (4,0) {$\cdots$};
    \draw[->] (4.3,0) to (4.8,0) ;
    \draw[->] (5.2,0) to node[midway, above]{\footnotesize{$p'_{-1}$}}(6.3,0) ;
    \draw[->, dotted] (6.7,0) to node[midway, above]{\footnotesize{$\beta'_2$}}(7.8,0) ;
    \draw[->] (1.7,-.8) to node[midway, above]{\footnotesize{$\alpha=w_1$}}(2.8,-.8) ;
    \draw[->] (3.2,-.8) to (3.7,-.8) ;
    \node at (4,-.8) {$\cdots$};
    \draw[->] (4.3,-.8) to (4.8,-.8) ;
    \draw[->] (5.2,-.8) to node[midway, below right=-.1]{\footnotesize{$w_m$}} (6.35,-.15) ;
\end{tikzpicture}\end{center}
and that $p_{-i}p_{-(i-1)},p'_{-i}p'_{-(i-1)}\in \cI$ for $1< i\leq n$, $w_{j}w_{j+1}\notin \cI$ for $1\leq i<m$. Moreover, if there exist arrows $\beta_i, \beta_i'$, assume that $\beta_1p_{-n}, p_{-1}\beta_2,\beta_1
'p_{-n}',p'_{-1}\beta_2'\notin\cI$, so that $p=p_{-n}\cdots p_{-1},p'=p'_{-n'}\cdots p'_{-1}\in\cL'$ and $w_1\cdots w_m=\gamma_r(\alpha)$. 
Then $(p,\epsilon_v)\in\Lplus$ and $(p', \gamma_r(\alpha))\in\Lpr{v'}$.
\end{example}

Note that the number of pairs $(p,\epsilon_v)$ in $\Lplus$ is at most $2-\outdeg(v)$, and the number of pairs $(p,\gamma_r(\alpha))$ is at most $\outdeg(v)$. Thus, $\abs{\Lplus}\leq 2$. If $(p,w), (p,w')\in \Lplus$, then $w'$ need not equal $w$ as seen in Example \ref{ex:rep-p}.

\begin{defn}
Let $I(v)$ be the $\ZZ$-graded right $A$-module with basis $\setst{p^{-1}}{p\in \cB, t(p)=v}$ where for $k\geq 0$,
\begin{equation*}
    (I(v))_{-k} = \kk\setst{p^{-1}}{p\in\cB_k, t(p)=v},
\end{equation*}
and for $p,q\in\cB$ with $t(p)=v$,
\begin{equation*}
    p^{-1}\cdot q = \begin{cases}
    r^{-1} & \text{if } p=qr \text{ for some }r\in\cB,\\
    0 & \text{otherwise}.
    \end{cases} 
\end{equation*}
Namely, for a path $p\in\cB_+$, $p^{-1}$ is its formal inverse. 
\end{defn}

Note that the submodule of $I(v)$ generated by $\epsilon_v$ is $S(v)$, and for a $k\geq 0$, an arbitrary element of $I(v)_{-k}$ has the form $\sum_{p\in A\epsilon_v\cap \cB_k}\lambda_pp^{-1}$ where $\lambda_p\in\kk$.

\begin{lemma}\label{lem:Iv-inj}
The module $I(v)$ is the graded injective hull of $S(v)$.
\end{lemma}

\begin{proof}
We begin by showing that $I(v)$ is an essential extension of $S(v)$. Let $M$ be a nonzero submodule of $I(v)$. Then, there exists a nonzero homogeneous element $x = \sum_{p\in A\epsilon_v\cap \cB_k}\lambda_pp^{-1}\in M_{-k}$ for some $k\geq 0,\lambda_p\in\kk$. Then, we have that $x\cdot \left(\sum_{p\in A\epsilon_v\cap \cB_k}\lambda_pp\right) = \left(\sum_{p\in A\epsilon_v\cap \cB_k}\lambda_p^2\right)\epsilon_v$, so $\epsilon_v\in S(v)\cap M$. It follows that $\epsilon_v$ is contained in any submodule of $I(v)$, and $I(v)$ is an essential extension of $S(v)$.  

Let $E$ be an essential extension of $I(v)$, and thus also of $S(v)$. Note that for $n>0$, since $S(v)_n=\set{0}$, if $z\in E_n$, $z$ generates an A-submodule of $E$, denoted $Z$, so that $S(v)\cap Z=\set{0}$. Therefore, it must be that $E_n=\set{0}$ for $n>0$. 

Assume for contradiction that $E\setminus I(v)\neq\set{0}$. Let $k$ be the minimal nonnegative integer $k\geq 0$ so that $(E\setminus I(v))_{-k}\neq \set{0}$. There must exist a vertex $v'$ (possibly $v$ itself) so that $(E\setminus I(v))_{-k}\epsilon_{v'}\neq \set{0}$.
Fix nonzero $z\in (E\setminus I(v))_{-k}\epsilon_{v'}$. If $k=0$, the submodule $z$ generates is $\kk z$, which has trivial intersection with $I(v)$, and thus we have a contradiction.

Otherwise, $k>0$. Then, by choice of $k$, $zA_1\subseteq E_{-k+1}\subseteq I(v)$ and thus
\begin{equation*}
    z\alpha=\sum_{p\in A\epsilon_v\cap \cB_{k-1}}\lambda_pp^{-1}.
\end{equation*}
Because $(z\alpha)p = \lambda_p\epsilon_v$ for any $p\in A\epsilon_v\cap \cB_{k-1}$,  $\lambda_p\neq 0$ implies that $\alpha p\notin \cI$. 
Thus, as $A$ is (locally) gentle, $z\alpha\neq 0$ implies that there is a unique  $p_\alpha\in A\epsilon_v\cap \cB_{k-1}$ for which $\lambda_{p_\alpha}\neq 0$. Consider 
\begin{equation*}
z'\defeq z-\sum_{\substack{\alpha\in\cQ_1\\z\alpha\neq 0}}\lambda_{p_\alpha} (\alpha p_\alpha)^{-1}.
\end{equation*}
By choice of of $p_\alpha$, we have $z'\in  (E\setminus I(v))_{-k}$ and $z'\alpha=0$ for every arrow $\alpha$. Moreover, as $z\in E\epsilon_{v'}$,  $z\alpha\neq 0$ implies that $s(\alpha)=v'$, and thus $z'\epsilon_{v''}=\delta_{v',v''}z'$. Thus the submodule of $E$ generated by $z'$ is $\kk z'$, which has trivial intersection with $I(v)$. This is a contradiction, hence $E=I(v)$.
\end{proof}

\begin{defn}\label{def:injtype2}
For an arrow $\alpha$, if $\gamma_\alpha$ is finite, let $I(\alpha)\defeq I\Big(t\big(\gamma_r(\alpha)\big)\Big)\big[-\ell(\gamma_r(\alpha))\big]$. If $\gamma_\alpha$ is infinite, let $I(\alpha)$ denote the $\ZZ$-graded right $A$-module with $\kk$-basis 
\begin{equation*}
    \setst{q^{-1}}{q\in\cB,L(q)=L(\gamma_r(\alpha))}\cup\{\epsilon_{s(\alpha)}\}\cup\setst{p\in\cB}{F(p)=\alpha},
\end{equation*}
i.e., initial subpaths and inverses of terminal subpaths of $\gamma_r(\alpha)^k$ for $k\geq 1$,  where 
\begin{equation*}
    ({I}(\alpha))_{-n}=\kk\setst{q^{-1}\in I(\alpha)}{\ell(q)=n}, \quad ({I}(\alpha))_n = \kk\setst{p\in I(\alpha)}{\ell(p)=n},
\end{equation*}
for $n>0$, $I_0=\kk\epsilon_{s(\alpha)}$, and the right module action of $A$ on ${I}(\alpha)$ is given by 
\begin{equation*}
    p\cdot w = pw,\quad
    q^{-1}\cdot w = \begin{cases}
    r^{-1} & \text{if }q=wr,\\
    r' & \text{if }w = qr',\\
    0 & \text{otherwise},
    \end{cases} \qquad 
    \epsilon_v\cdot w= \begin{cases}
    w & \text{if }F(w)=\alpha\text{ or } w=\epsilon_v,\\ 
    0 & \text{ otherwise,}
    \end{cases} 
\end{equation*}
for a path $p\in I(\alpha)_+$, an inverse path $q^{-1}\in I(\alpha)_{<0}$, $r\in \cB_+$ and $w, r'\in \cB$. 
\end{defn}

Note that the action of $A$ on $q^{-1}\in\cI(\alpha)$ depends on whether or not $\gamma_\alpha$ is infinite. In particular, if there exists $\beta$ such that $q\beta\in\cB$, then $q^{-1}\cdot (q\beta)$ is zero if $\gamma_\alpha$ is finite, but
$\beta$ if $\gamma_\alpha$ is infinite.

If $\gamma_\alpha$ is infinite, then $\dim_\kk\left(I(\alpha)_n\right)=1$ for every $n\in \zz$. If $\alpha$ and $\beta$ are contained in the same  maximal path, then $I(\alpha)$ and $I(\beta)$ are isomorphic, up to a shift. In particular, if $\gamma_r(\alpha)=\alpha_1\ldots\alpha_n$, then $I(\alpha_i)\cong I(\alpha)[i-1]$
for $1\leq i\leq n$,  with isomorphism given by $x\mapsto \alpha_1\dotsc \alpha_{i-1}x$ where this left multiplication is the natural one.

\begin{lemma} 
The module $I(\alpha)$ is the graded injective hull of $\alpha A$. 
\end{lemma}

\begin{proof}
The case where $\gamma_\alpha$ is finite follows from Lemma \ref{lem:Iv-inj}, as $\kk\gamma_r(\alpha)\cong S(t(\gamma_r(\alpha)))[\ell(\gamma_r(\alpha))]$ is an essential submodule of $\alpha A$.
Suppose that $\gamma_\alpha$ is infinite, and let $\gamma_r(\alpha)=\alpha_1\dotsc \alpha_\ell$ and $\alpha_{i+j\ell}\defeq \alpha_i$ for any $1\leq i\leq n$ and $j\in\ZZ$. Let 
\begin{equation*}
    \baselt_n\defeq\begin{cases} \alpha_{\ell}^{-1}\alpha_{\ell-1}^{-1}\dotsc \alpha_{\ell+n+1}^{-1}&\text{if }n<0,\\ \epsilon_{s(\alpha)} & \text{if }n=0,\\
    \alpha_1\alpha_2\dotsc\alpha_{n} & \text{if }n>0,\end{cases}
\end{equation*}
so that $(I(\alpha))_n=\kk \set{\baselt_n}$ and $\baselt_n\alpha_{n+1}=\baselt_{n+1}$ for any $n\in\ZZ$. Then, for any $z\in I(\alpha)\setminus\set{0}$ there exist $\lambda_i\in\kk$ so that $z=\sum_i \lambda_i\baselt_i$ and a minimal $n\in \ZZ$ for which $\lambda_n\neq 0$. Either $n>0$ and  $z\in \alpha A$, or $z\cdot \left(\alpha_{\ell+n+1}\cdots \alpha_{\ell-1}\alpha_\ell\alpha_1\right)\in \lambda_n\alpha_1+\alpha A_+$. Thus $\alpha A$ is an essential submodule of $I(\alpha)$.
    
Suppose that $E$ is an essential extension of $I(\alpha)$ and assume for contradiction that there exists $z\in E_n\setminus I(\alpha)$. 
Then there exists $a\in A_k$ for some $k\geq 0$ such that $z\cdot a\in I(\alpha)_{n+k}\setminus\{0\}$. Then $z\cdot a$ is a nonzero multiple of $\baselt_{n+k}$, so
$(z\cdot a)\cdot \alpha_{n+k+1}\neq 0$ 
and thus $a\cdot \alpha_{n+k+1}\notin \cI$. 
As $a\in A_k$, this implies that $a\cdot \alpha_{n+k+1}$ is a nonzero multiple of $\alpha_{n+1}\dotsc\alpha_{n+k+1}$, so $z\cdot (\alpha_{n+1}\dotsc\alpha_{n+k+1})\in I(\alpha)\setminus\{0\}$ as well. By scaling $z$, we may assume without loss of generality that 
\begin{equation*}
    z\cdot (\alpha_{n+1}\dotsc\alpha_{n+k+1})=\baselt_{n+k+1}.
\end{equation*}
In particular, whenever $m>k$,
\begin{equation*}
    z\cdot (\alpha_{n+1}\dotsc\alpha_{n+m})
    =\baselt_{n+m}=\baselt_n\cdot (\alpha_{n+1}\dotsc\alpha_{n+m}).
\end{equation*}
Since $z\notin I(\alpha)$, we have $z-\baselt_n\in E\setminus\{0\}$. Hence there exists $a'\in A_{k'}$ such that $(z-\baselt_n)\cdot a'\in I(\alpha)\setminus\{0\}$. Then
$\baselt_n\cdot a'=\gamma \baselt_{n+k'}$,
where $\gamma$ is the coefficient of $\alpha_{n+1}\dotsc\alpha_{n+k'}$ in $a'$.
As $\baselt_n\cdot a'\in I(\alpha)$, we have $z\cdot a'\in I(\alpha)$ as well, so 
$z\cdot a'=\gamma'\baselt_{n+k'}$
for some $\gamma'\in \kk$.
Take $m=\max\{k,k'\}+1$. Then
\begin{equation*}
    (z\cdot a')\cdot (\alpha_{n+k'+1} \dotsc \alpha_{n+m})=\gamma'\baselt_{n+k'}\cdot(\alpha_{n+k'+1}\dotsc\alpha_{n+m})=\gamma'\baselt_{n+m}.
\end{equation*}
As $a'\alpha_{n+k'+1}\dotsc \alpha_{n+m}=\gamma(\alpha_{n+1}\dotsc\alpha_{n+m})$ and $m>k$, 
\begin{align*} 
    (z\cdot a')\cdot (\alpha_{n+k'+1} \dotsc \alpha_{n+m})=  z\cdot (a'\alpha_{n+k'+1}\dotsc \alpha_{n+m})=\gamma z\cdot (\alpha_{n+1}\dotsc \alpha_{n+m})
    =\gamma \baselt_{n+m}
\end{align*}
as well, so $\gamma=\gamma'$. Then $(z-\baselt_n)\cdot a'=(\gamma'-\gamma)\baselt_{n+k'}=0$, a contradiction.  Therefore $E=I(\alpha)$.
\end{proof}

\begin{lemma}\label{lem:injproj} If $\outdeg(v)=1$ and $\Lplus=\emptyset$, then $\epsilon_vA\cong I(\alpha)$, where $\alpha$ is the unique arrow with source $v$.\end{lemma}
\begin{proof}
Assume that $\alpha$ is the only arrow with source $v$. 
If $\gamma_r(\alpha)$ is not left maximal, then there exists an arrow $\beta$ such that $\beta\alpha\notin\cI$ and $(\gamma_\beta^\#,\epsilon_v)\in \Lplus$. Recall that if  $\gamma_r(\alpha)$ is left maximal, then $\gamma_\alpha=\gamma_r(\alpha)$. If $\gamma_r(\alpha)$ is right maximal and $\outdeg(t(\gamma_r(\alpha)))=2$, then
there exists $\beta\neq L(\gamma_r(\alpha))$ such that $t(\beta)=t(\gamma_r(\alpha))$ and
$(\gamma_\beta^\#,\gamma_r(\alpha))\in \Lplus$. Consequently, $\Lplus=\emptyset$ implies that $\gamma_\alpha=\gamma_r(\alpha)$ and $\indeg(t(\gamma_\alpha))=1$. 

Assume that $\Lplus=\emptyset$. Then 
\begin{equation*}
I(\alpha)=I(t(\gamma_\alpha))[-\ell(\gamma_\alpha)]=\kk\setst{q^{-1}}{\gamma_\alpha=rq\text{ for some }r\in\cB}
\end{equation*}
as $t(q)=t(\gamma_\alpha)$ implies that $L(q)=L(\gamma_\alpha)$ and thus $q$ is a terminal subpath of $\gamma_\alpha$. Consider the graded morphism $\phi: I(\alpha)\rightarrow \epsilon_vA$ given by $\phi(q^{-1})=r$ for paths $r,q\in\cB$ such that $\gamma_\alpha=rq$. This $\phi$ is an isomorphism, as $\outdeg(v)=1$ implies that 
\begin{equation*} \epsilon_vA=\kk\setst{r}{r\text{ is an initial subpath of $\gamma_r(\alpha)^k$ for some }k\geq 0}\end{equation*}
and $\gamma_\alpha=\gamma_r(\alpha)$.
\end{proof}

\begin{defn}
Suppose that $w\in \cB_+$ and $k\in\ZZ$. Define $w^*: I(t(w))[k]\rightarrow I(s(w))[k+\ell(w)]$ to be the homomorphism of graded right $A$-modules given by
\begin{equation*}
    w^*(q^{-1})=
    \begin{cases}
        r^{-1} & \text{ if }q=rw,\\
        0 & \text{ otherwise}.
    \end{cases}
\end{equation*}
\end{defn}
The image of the homomorphism $w^*$ is the span of $q^{-1}$ such that $qw\in\cB$ and its kernel is the span of $q^{-1}$ such that $w$ is not a terminal subpath of $q$ but $t(w)=t(q)$. That is, \begin{align*}\im (w^*)&=\kk\setst{q^{-1}}{q\in\cB, qw\notin\cI},\\
\Ker(w^*)&=\kk\setst{q^{-1}}{q\in\cB, t(q)=t(w), q\neq q'w\text{ for any }q'\in \cB}
.\end{align*}

For the graded injective resolution of $A\epsilon_v$, we will use the notation of the following lemma.
\begin{lemma}\label{lem:jpw}
    For $(p,w)\in \Lplus$, let 
\begin{equation} I^j(p,w)\defeq \begin{cases}
I(t(p_{-1}))[-\ell(w)] &\text{ if } j=0,\\
I(s(p_{-j}))[j-\ell(w)] & \text{ if }0<j\leq \ell(p),\\
\end{cases}\end{equation}
and, for $0\leq j<\ell(p)$, let
\begin{equation} d^j_{(p,w)}\defeq p^*_{-(j+1)}:I^j(p,w)\rightarrow I^{j+1}(p,w).\end{equation}
Then 
\begin{tikzcd}
\Ker(d^0_{(p,w)})\hookrightarrow {I^0(p,w)} \arrow[r, "{d^0_{(p,w)}}"] & {I^1(p,w)} \arrow[r, "{d^1_{(p,w)}}"] & \cdots \arrow[r, "{d^{\ell(p)-1}_{(p,w)}}"] & {I^{\ell(p)}(p,w)}       \arrow[r, ""]      &0
\end{tikzcd}
is  exact. 
\end{lemma}
\begin{proof}
Suppose that $q\in\cB$ has target $s(p_{-j})$ for some $j$, so that $q^{-1}\in I^j(p,w)$. As $p_{-(j+1)}p_{-j}\in\cI$ and $A$ is (locally) gentle, $q$ ends with the arrow $p_{-(j+1)}$ if and only if $qp_{-j}\in\cI$.
 Thus, for $1\leq j<\ell(p)$, we have 
 \begin{equation*}
     \im (d^{j-1}_{(p,w)})=\kk\setst{q^{-1}}{ qp_{-j}\notin\cI}=\kk\setst{q^{-1}}{t(q)=t(p_{-(j+1)}), L(q)\neq p_{-(j+1)}} = \Ker(d^{j}_{(p,w)}).
 \end{equation*} 
As $p$ is maximal, we know that $qp_{-\ell(p)}\notin\cI$ for any path $q\in\cB$ with target $s(p_{-\ell(p)})$, so
 \begin{equation*}
     \im (d^{\ell(p)-1}_{(p,w)})=\kk\setst{q^{-1}} {qp_{-\ell(p)}\notin\cI}=\kk\setst{q^{-1}}{q\in\cB, t(q)=s(p_{-\ell(p)})}=I^{\ell(p)}(p,w).
     \end{equation*}

\end{proof}

\begin{lemma}\label{lem:injres}
Let $A=\kk\cQ/\cI$ be a (locally) gentle algebra and $v\in \cQ_0$. 
Then there is a graded injective resolution of length $m$ for $\epsilon_vA$, where
\begin{equation*}
    m=\begin{cases}
    \max\setst{\ell(p)}{ (p,w)\in \Lplus} & \text{if }\Lplus\neq \emptyset,\\
    1 & \text{if }\Lplus=\emptyset,\,\outdeg(v)=2,\\
    0 & \text{otherwise}.
    \end{cases}
\end{equation*}
\end{lemma}
In particular,  the module $\epsilon_v A$ has a graded injective resolution of length $m$ for which the first term is determined by $s^{-1}(v)$ and the higher terms trace out forbidden paths $p$ where $(p,w)\in\Lplus$. This resolution is defined as follows.
\begin{defn} \label{def:injres} Let $m$ be as in Lemma \ref{lem:injres}. Consider the resolution of $\epsilon_vA$ given by
\begin{equation*}
    0\rightarrow \epsilon_vA\xrightarrow{\iota} I^0\xrightarrow{d^0} I^1 \xrightarrow{d^1}\dotsc \xrightarrow{d^{m-1}} I^m\rightarrow 0,
\end{equation*}
where the modules and maps are defined as follows. 
The modules $I^j$ are given by
\begin{equation*}
    I^0\defeq  \begin{cases}
    I(v) & \text{if }\outdeg(v)=0,\\
    \bigoplus\limits_{\alpha\in s^{-1}(v)} I(\alpha) & \text{if }\outdeg(v)\neq 0,
    \end{cases}
\end{equation*}
and for $j>0$,
\begin{equation*}
    I^j\defeq \begin{cases} 
    I(v)\oplus \bigoplus\limits_{(p,w)\in \Lplus} I^j(p,w) & \text{ if } j=1, \outdeg(v)=2,\\
    \bigoplus\limits_{\substack{(p,w)\in \Lplus\\\ell(p)\geq j}} I^j(p,w)
    &\text{ otherwise}
    \end{cases}
\end{equation*}
where  $I^j(p,w)$ are as in Lemma \ref{lem:jpw}.
In particular, $I^j=0$ for $j>m$. If $v$ is a sink, then $\epsilon_vA\cong S(v)$ and $\iota$ is the natural inclusion. Otherwise, 
\begin{equation*}
\iota = \begin{cases}(\iota_\alpha) & \text{ if }s^{-1}(v)=\set{\alpha},\\
\left(\begin{matrix}\iota_\alpha&\iota_\beta\end{matrix}\right)^T & \text{ if }s^{-1}(v)=\set{\alpha,\beta},
\end{cases}
\end{equation*}
where, for $\alpha\in s^{-1}(v)$, $\iota_\alpha: \epsilon_vA\rightarrow I(\alpha)$ is given by 
\begin{equation*} 
    \epsilon_v\mapsto
    \begin{cases} \gamma_r(\alpha)^{-1}&\text{ if }\ell(\gamma_\alpha)<\infty\\
    \epsilon_v&\text{ if }\ell(\gamma_\alpha)=\infty.\end{cases}
\end{equation*}
For $\alpha\in s^{-1}(v)$, let $\rho_\alpha: I(\alpha)\rightarrow I(v)$ be the composition $I(\alpha)\xrightarrow{\pi}I(\alpha)/(I(\alpha))_+\hookrightarrow I(v)$ if $\gamma_\alpha$ is infinite, or $\rho_\alpha=(\gamma_r(\alpha))^*$ if $\gamma_\alpha$ is finite.
If $\outdeg(v)=2$, let $\rho\defeq \left(\begin{matrix} \rho_{\alpha_1} & -\rho_{\alpha_2}\end{matrix}\right)$ where  $s^{-1}(v)=\set{\alpha_1,\alpha_2}$. 

Let
\begin{equation*}
    d^0\defeq \begin{cases} 
   \left(\begin{matrix}\rho &\bigoplus\limits_{(p,w)\in \Lplus}d^0_{(p,w)}\end{matrix}\right)^T
    & \text{ if } \outdeg(v)=2,\\
   \left(\begin{matrix} \hat{d}^0_{(p,w)}&\hat{d}^0_{(p',w')}\end{matrix}\right)^T&  \text{ if }\Lplus=\set{(p,w), (p',w')},\\
    (\hat{d}^0_{(p,w)})& \text{ if } \Lplus=\set{(p,w)},\\
    \end{cases}
\end{equation*}
 where $
    \hat{d}^0_{(p,w)} \defeq 
    \begin{cases} 
    d^0_{(p,w)}\circ \rho_\alpha & \text{ if  $w=\epsilon_v$,  $ s^{-1}(v)=\set{\alpha}$,}\\
    d^0_{(p,w)}&\text{ otherwise,}
    \end{cases}$ and, for $i\geq 1$, let 
\begin{equation*}
    d^i\defeq \bigoplus\limits_{\substack{(p,w)\in \Lplus\\ \ell(p)>i}}d^i_{(p,w)}.
\end{equation*}
\end{defn}

Before proving that this is a graded injective resolution, we provide some examples to better understand this resolution.

\begin{example}
    Consider the gentle algebra $A=\kk\cQ/\langle c_1\alpha_1, \alpha_1\alpha_2, \alpha_2 c_2\rangle$, where $\cQ$ is the quiver
\begin{center}
\begin{tikzpicture}
    \vtx{0,0}{1}
    \vtx{2,0}{2}
    \vtx{4,0}{3}
    \draw[->] (.2,0) to node[midway, above]{\footnotesize{$\alpha_1$}}(1.8,0) ;
    \draw[->] (2.2,0) to node[midway, above]{\footnotesize{$\alpha_2$}}(3.8,0) ;
    \draw[->] (-.1,.15) to[out=120, in=60, looseness=10] node[midway, above]{\footnotesize{$c_1$}}
    (.1,.15) ;
    \draw[->] (3.9,.15) to[out=120, in=60, looseness=10] node[midway, above]{\footnotesize{$c_2$}}
    (4.1,.15) ;
\end{tikzpicture}.
\end{center}
Here $\cL'=\{c_1\alpha_1\alpha_2 c_2\}$, and $p\defeq c_1\alpha_1\alpha_2 c_2$ has target 3. The only arrows with finite right maximal paths are $\alpha_i$, and $\gamma_r(\alpha_i)=\alpha_i$. 
As the target of the one finite maximal forbidden path $p$ is not $1$ or $t(\alpha_1)$, we have $\Lpr{1}=\emptyset$. The other pair sets are
$\Lpr{2}=\{(p, \alpha_2)\}$ and $\Lpr{3}=\{(p, \epsilon_3)\}$. The resolutions described in Definition \ref{def:injres} are as follows:
\begin{align*}
    &0\rightarrow \epsilon_1A\xrightarrow{\epsilon_1\mapsto  \left(\begin{matrix}\scriptstyle(\alpha_1)^{-1}\\\scriptstyle\epsilon_1\end{matrix}\right)} \underbrace{I(2)[-1]}_{I(\alpha_1)}\oplus I(c_1)\xrightarrow{\rho = \left(\begin{matrix}
       \scriptstyle \alpha_1^*&\scriptstyle-\rho_{c_1} 
    \end{matrix}\right)\cdot -} I(1)\rightarrow 0,\\
    &0\rightarrow \epsilon_2A\xrightarrow{\epsilon_2\mapsto (\alpha_2)^{-1}} \underbrace{I(3)[-1]}_{I(\alpha_2)}
    \xrightarrow{c_2^*}  I(3)\xrightarrow{\alpha_2^*}I(2)[1]\xrightarrow{\alpha_1^*}I(1)[2]\xrightarrow{c_1^*}I(1)[3]\rightarrow 0,\\
    &0\rightarrow \epsilon_3A\xrightarrow{\epsilon_3\mapsto \epsilon_3} I(c_2)\xrightarrow{(c_1^*)\circ {\displaystyle\rho_{c_2}}}  I(3)[1]\xrightarrow{\alpha_2^*}I(2)[2]\xrightarrow{\alpha_1^*}I(1)[3]\xrightarrow{c_1^*}I(1)[4]\rightarrow 0.\\
\end{align*}
The resolutions  $\epsilon_2A$ and $\epsilon_3A$ both trace out the same maximal forbidden thread, $p$, but the shift differs as the corresponding $w$ differs. Additionally, in the resolution of $\epsilon_3 A$, we have $d^0={d}^0_{(p,\epsilon_3)}\circ \rho_{c_2}$ as $c_1$ is the only arrow with source 3 and the corresponding $w$ is not $\gamma_r(c_1)$.
\end{example}

In the following example, we illustrate some cases in which the same forbidden path occurs in two pairs in $\Lplus$.
\begin{example}\label{ex:rep-p}
Consider the gentle algebra $A=\kk\cQ/\langle \beta_2\beta_1\rangle$, where $\cQ$ is the quiver
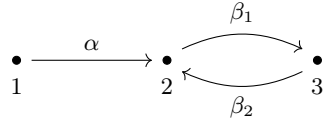
\begin{center}
\begin{tikzpicture}
    \vtx{0,0}{1}
    \vtx{2,0}{2}
    \vtx{4,0}{3}
    \draw[->] (.2,0) to node[midway, above]{\footnotesize{$\alpha$}}(1.8,0) ;
    \draw[->] (2.2,.15) to[out =25, in=155] node[midway, above]{\footnotesize{$\beta_1$}}(3.8,.15) ;
    \draw[<-] (2.2,-.15) to[out =-25, in=205] node[midway, below]{\footnotesize{$\beta_2$}}(3.8,-.15) ;
\end{tikzpicture}.
\end{center}
Then $\Lpr{2}=\{(\alpha,\epsilon_2), (\alpha, \beta_1\beta_2)\}$ and the resolution for $\epsilon_2A$ described in Definition \ref{def:injres} is 
\begin{equation*}
    0\rightarrow \epsilon_2A\xrightarrow{\epsilon_2\mapsto (\beta_1\beta_2)^{-1}} I(2)[-2]\xrightarrow{\left(\begin{matrix}\scriptstyle(\alpha\beta_1\beta_2)^*\\ \scriptstyle\alpha^*\end{matrix}\right)\cdot -}  I(1)[1]\oplus I(1)[-1]\rightarrow 0.
\end{equation*}

Consider $A=\kk\cQ/\langle c^2\rangle$ where $\cQ$ is the quiver
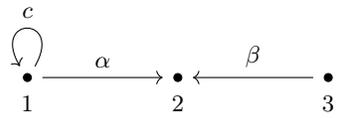
\begin{center}
\begin{tikzpicture}
    \vtx{0,0}{1}
    \vtx{2,0}{2}
    \vtx{4,0}{3}
    \draw[->] (.2,0) to node[midway, above]{\footnotesize{$\alpha$}}(1.8,0) ;
    \draw[<-] (2.2,0) to node[midway, above]{\footnotesize{$\beta$}}(3.8,0) ;
    \draw[<-] (-.1,.15) to[out=120, in=60, looseness=10] node[midway, above]{\footnotesize{$c$}}
    (.1,.15) ;
\end{tikzpicture}.
\end{center}
Then $\Lpr{1} =\{(\beta,\alpha),(\beta,c\alpha)\}$ and the resolution for $\epsilon_1A$ described in Definition \ref{def:injres} is 
\begin{equation*}
    0\rightarrow \epsilon_1A\xrightarrow{\epsilon_1\mapsto\left(\begin{matrix} \scriptstyle \alpha^{-1}\\\scriptstyle (c\alpha)^{-1}\end{matrix}\right)} I(2)[-1]\oplus I(2)[-2] \xrightarrow{\left(\begin{matrix}\scriptstyle\alpha^* & \scriptstyle-(c\alpha)^*\\ \scriptstyle\beta^* & \scriptstyle 0\\ \scriptstyle0 & \scriptstyle\beta^*\end{matrix}\right)\cdot -} I(1)\oplus I(1)[0]\oplus I(1)[-1]\rightarrow 0.
\end{equation*}
\end{example}

We will now show that the resolution of Definition \ref{def:injres} is indeed a graded injective resolution of $\epsilon_vA$, proving Lemma \ref{lem:injres}.

\begin{proof} We have already shown the case $m=0$ in Lemma \ref{lem:injproj}, noting that $\phi^{-1}=\iota_\alpha$ in that case. By Lemma \ref{lem:jpw}, it only remains to show that $0\rightarrow \epsilon_vA\xrightarrow{\iota} I^0\xrightarrow{d^0}I^1\xrightarrow{d^1}\im d^1$ is exact.

\textbf{Case 1 ($\outdeg(v)=0$):} We have $\epsilon_vA\cong S(v)$ and $\iota$ is the usual inclusion, so $\im(\iota)=\kk\epsilon_v$. As $v$ is a sink, $\gamma^\#_l(\beta)=\gamma_\beta^\#$ for all arrow $\beta$ with target $v$ and thus $\Lplus(v)=\setst{(\gamma_\beta^\#,\epsilon_v)}{t(\beta)=v}$. As
$\Ker d^0_{(\gamma_\beta^\#,\epsilon_v)}=\kk\setst{q^{-1}}{q\in\cB, t(q)=v, L(q)\neq \beta}$, we see that
\begin{equation*}\Ker d^0=\bigcap_{t(\beta)=v}\Ker d^0_{(\gamma_\beta^\#,\epsilon_v)}=\kk\epsilon_v=\im(\iota).
\end{equation*}
If $t^{-1}(v)=\{\beta_1,\beta_2\}$, then $q\notin  \Ker d^0_{(\gamma_{\beta_1}^\#,\epsilon_v)}$ implies $q\in  \Ker d^0_{(\gamma_{\beta_2}^\#,\epsilon_v)}$. Thus
\begin{equation*}
    \im d^0=\bigoplus_{t(\beta)=v} \im(d^0_{(\gamma_{\beta}^\#,\epsilon_v)}),
\end{equation*}
which is equal to $\Ker d^1$ by Lemma \ref{lem:jpw}.

\textbf{Case 2 ($\outdeg(v)=1$):} If $\alpha$ is the only arrow with source $v$, then $\iota=\iota_\alpha$ is injective and
\begin{equation*}
    \im(\iota_\alpha)=\begin{cases}
    \kk\setst{q^{-1}}{\text{q is  a terminal subpath of $\gamma_r(\alpha)$}} & \text{if }\ell(\gamma_\alpha)<\infty,\\
    I(\alpha)_{\geq 0} & \text{if }\ell(\gamma_\alpha)=\infty.
    \end{cases}
\end{equation*}
If $(p,\epsilon_v)\in\Lplus$, then $p_{-1}\alpha\notin \cI$ and thus 
\begin{equation*}
\im(\rho_\alpha)=\kk\setst{q^{-1}}{q\alpha\notin\cI}=\kk\left(\setst{q^{-1}}{q\in\cB, L(q)=p_{-1}}\cup\set{\epsilon_v^{-1}}\right).
\end{equation*}
As
$\Ker(d^0_{(p,\epsilon_v)})=\kk\setst{q^{-1}}{q\in\cB, t(q)=v, L(q)\neq p_{-1}}$, we see that $\Ker(\hat{d}^0_{(p,\epsilon_v)})=\Ker(\rho_\alpha)\cup\im(\kk\epsilon_v)$ and $\im (\hat{d}^0_{(p,\epsilon_v)})=\im ({d}^0_{(p,\epsilon_v)})$. The latter implies that $\im(d^0)=\bigoplus_{(p,w)\in\Lplus}\im(d^0_{(p,w)})=\Ker(d^1)$ by Lemma \ref{lem:jpw}.

Note that
\begin{equation*}
    \iota_{\alpha}(\epsilon_v A_+)\subseteq \Ker(\rho_{\alpha})=\begin{cases} \kk\setst{q^{-1}}{\text{$\gamma_r({\alpha})$ is not a terminal subpath of $q$}} & \text{if $\ell(\gamma_r(\alpha))<\infty$},\\
    I(\alpha)_+ & \text{if $\ell(\gamma_r(\alpha)) =\infty$}.
    \end{cases}
\end{equation*}
If $\iota_{\alpha} (\epsilon_vA_+) \neq \Ker(\rho_{\alpha})$, then $\gamma_{\alpha}$ is finite and there exists an arrow  $\beta$ such that $(\gamma^\#_{\beta},\gamma_r(\alpha))\in \Lplus$. In that case, $\Ker (d^0_{(\gamma^\#_{\beta},\gamma_r(\alpha))})=\kk\setst{q^{-1}}{q\in\cB,  L(q)=L(\gamma_r(\alpha))}$, so 
\begin{align*}
\Ker(d^0)&=\Ker(\hat{d}^0_{(p,\epsilon_v)})\cap \Ker(d^0_{(\gamma^\#_{\beta},\gamma_r(\alpha))})\\&=    \kk\setst{q^{-1}}{\text{q is  a terminal subpath of $\gamma_r(\alpha)$}} =\im (\iota).
\end{align*}

If $\Lplus\neq \emptyset$ but $(p,\epsilon_v)\notin\Lplus$ for any $p\in\cL'$, then there exists $p\in\cL'$ so that $\Lplus=\set{(p,\gamma_r(\alpha))}$. As there is no maximal forbidden path with target $v$, $\gamma_r(\alpha)$ is right maximal and thus
\begin{align*}
\Ker(d^0)=\Ker(d^0_{(p,\gamma_r(\alpha))})&=\kk\setst{q^{-1}}{q\in\cB,  L(q)=L(\gamma_r(\alpha))}\\&=\kk\setst{q^{-1}}{\text{q is  a terminal subpath of $\gamma_r(\alpha)$}}=\im(\iota_\alpha).
\end{align*}

\textbf{Case 3 ($\outdeg(v)=2$):} Suppose $\alpha_1\neq \alpha_2$ both have source $v$. As \begin{equation*}
    \Ker(\iota) = \Ker(\iota_{\alpha_1})\cap \Ker(\iota_{\alpha_2})=\alpha_2 A\cap \alpha_1 A=0,
\end{equation*}
 $\iota$ is again injective. Since $\outdeg(v)=2$, $(p,w)\in \Lplus$ implies that $w=\gamma_r(\alpha_i)$ for some $i$. If $(p,\gamma_r(\alpha_i))\in \Lplus$, then $d^0_{(p,\gamma_r(\alpha_i))}\circ \iota_{\alpha_i} =0$ by our assumption that $p_{-1}\neq L(\gamma_r({\alpha_i}))$.
Since $\im(\rho_{\alpha_i})=\kk\setst{q^{-1}}{q{\alpha_i}\notin\cI}$, we see that $\im(\rho_{\alpha_1})\cap \im(\rho_{\alpha_2})=\kk\epsilon_v^{-1}=(\rho_{\alpha_i}\circ\iota_{\alpha_i})(\epsilon_v)$ and thus
\begin{equation*}
    \Ker(\rho) = \setst{(x,y)\in I(\alpha_1)\oplus I(\alpha_2)}{\rho_{\alpha_1}(x) =\rho_{\alpha_2}(y)}=\kk\iota(\epsilon_v)+\left(\Ker(\rho_{\alpha_1})\oplus \Ker(\rho_{\alpha_2})\right).
\end{equation*}
 Since $\iota_{\alpha_i}(\epsilon_v A_+)\subseteq \Ker(\rho_{\alpha_i})$, we know that  $\rho\circ \iota=0$. As in the previous case, if 
 $\iota_{\alpha_i} (\epsilon_vA_+) \neq \Ker(\rho_{\alpha_i})$, then there exists an arrow $\beta_i$ such that $(\gamma^\#_{{\alpha_i}},\gamma_r(\alpha_i))\in \Lplus$ and  
 $\Ker(\rho_{\alpha_i})\cap \Ker(d^0_{(\gamma^\#_{\beta_i},\gamma_r(\alpha_i))})=\iota_{\alpha_i}(\epsilon_vA_+)$. Thus $\Ker (d^0)=\Ker (\rho)\cap\bigcap_{(p,w)\in \Lplus}\Ker (d^0_{(p,w)})=\im (\iota)$.

Note that $\im(\rho)=I(v)$, since $t(q)=v$ implies that $q\alpha_i$ does not belong to $\cI$ for at least one $i$. Thus $\im d^0=\Ker d^1$ by Lemma \ref{lem:jpw}. 
\end{proof}

\begin{example}\label{ex:injresforASG} In this example, we compute the injective resolutions of $A$ in two special cases. 
Suppose $A=\kk\widetilde{A}_n$. For $i\in \cQ_0$, we have $I(\alpha_i)=I(i+1)[-1]$ and $\Lpr{i}=\{(\alpha_{i-1},\epsilon_i)\}$, so
\begin{equation*}
    \epsilon_iA\xrightarrow{[\alpha_i^{-1}]} I(i+1)[-1]\xrightarrow{(\alpha_{i-1}\alpha_i)^*}I(i-1)[1] \rightarrow 0 
\end{equation*}
is an injective resolution of $\epsilon_iA$ by Lemma \ref{lem:injres} and thus
\begin{equation*}
    A_A\hookrightarrow  \bigoplus_{0\leq i\leq n} I(i+1)[-1]\rightarrow \bigoplus_{0\leq i\leq n} I(i-1)[1]\rightarrow 0  
\end{equation*}
is an injective resolution of $A$.
    
Suppose instead that $\cQ$ is 2-regular. Then $A^\#$ has no finite maximal paths, so $\Lplus=\emptyset$ for every vertex $v$. If $\alpha$ and $\beta$ are the arrows with source $v$, then 
\begin{equation*}
    \epsilon_vA\hookrightarrow I(\alpha)\oplus I(\beta)\rightarrow I(v)\rightarrow 0 
\end{equation*}
is an injective resolution of  $\epsilon_vA$ by Lemma \ref{lem:injres} and thus
\begin{equation*}
    A_A\hookrightarrow  \bigoplus_{\alpha\in\cQ_1} I(\alpha)\rightarrow  \bigoplus_{v\in\cQ_0} I(v)\rightarrow 0  
\end{equation*}
is an injective resolution of $A$.
\end{example}

The following theorem extends  \cite[Theorem 6.9]{LGH} to locally gentle algebras.
\begin{theorem}\label{thm:injdim}
Let $A=\kk\cQ/\cI$ be (locally) gentle. 
Then the graded injective dimension of $A$ is
\begin{equation*}
    \injdim(A)=\begin{cases}
    \max\setst{\ell(p)}{p\in \cL'} & \text{if }\cL'\neq \emptyset,\\
    0 & \text{if }A=(\kk\widetilde{A}_n)^\#,\\
    1 & \text{otherwise}.
    \end{cases}
\end{equation*}
In particular, the graded injective resolutions given in Example \ref{ex:injresforASG} are minimal.
\end{theorem}

\begin{proof}
Suppose that $\cL'\neq \emptyset$ and let $M=\max\{\ell(p):p\in \cL'\}$. Note that $M\geq 1$. For any $v\in\cQ_0$, if $(p,w)\in \Lplus$, then $p\in\cL'$, and thus, $\ell(p)\leq M$, so $\injdim(\epsilon_vA)\leq M$ by Lemma \ref{lem:injres}. Thus $\injdim(A_A)=\max\setst{\injdim(\epsilon_vA)}{v\in\cQ_0}\leq M$. There exists $p\in \cL'$ such that $M=\ell(p)$. Then $\grExt^M_A(S(s(p)),A)\neq 0$ by Corollary \ref{cor:extcalc} and so $\injdim(A_A)\geq M$. Therefore, $\injdim(A_A)= M$.

Now, suppose $\cL'=\emptyset$. Then $\Lplus=\emptyset$ for each vertex $v$ as well, so $\injdim(A)\leq 1$ by Lemma \ref{lem:injres}.
As $A^\#$ has no finite maximal paths,  either $A^\#=\kk \widetilde{A}_n$ or $A$ has a vertex which is the source of two arrows by Lemma \ref{lem:invoutdeg}. Suppose that $A=(\kk\widetilde{A}_n)^\#$. Then $\injdim(\epsilon_vA)=0$ for each vertex $v$ by Lemma \ref{lem:injres} and thus $\injdim(A)=0$. Otherwise, let $v$ be a vertex which is the source of two arrows.  Then $\grExt^1_A(S(v),A)\neq 0$ and thus $\injdim(A_A)\geq 1$ by Corollary \ref{cor:extcalc}, whence $\injdim(A_A)=1$. 
    
As the opposite algebra, $A^\text{op}$, is also a (locally) gentle algebra, 
$\injdim(_AA)=\injdim(A^\text{op}_{A^\text{op}})<\infty$. As $A$ is Noetherian, $\injdim(_AA),\injdim(A_{A})<\infty$ implies $\injdim(_AA)=\injdim(A_A)$ \cite{Zaks}.
\end{proof}

\section{Homological conditions}\label{sec:ASconds}
\subsection{Artin-Schelter Gorenstein and Regular}

Recall that $S=A/\grJ(a)\cong \bigoplus_{v\in\cQ_0}S(v)$.
In this section, we use the results from the previous section to determine AS conditions for locally gentle algebras.  We say that $\gA$ is \emph{Artin-Schelter Gorenstein} (AS Gorenstein) of dimension $d$ if $\injdim(_\gA\gA)=\injdim(\gA_\gA)=d<\infty$ and
\begin{equation*}
\grExt^i_\gA(S,\gA)\cong 
    \begin{cases}
        V[\ell] & \text{if }i=d,\\
        0 & \text{otherwise},
    \end{cases}
\end{equation*}
where $V$ is a $\kk$-central invertible $(S,S)$-bimodule concentrated in degree 0. The value $\ell$ is referred to as the \emph{Gorenstein parameter} of $\gA$. We say that $\gA$ is \emph{Artin-Schelter  regular} (AS regular) of dimension $d$ if $\lgldim(\gA)=\rgldim(\gA)=d<\infty$) and
\begin{equation*}
\grExt^i_\gA(S,\gA)\cong \begin{cases}
    V[\ell] & \text{if }i=d,\\
    0 & \text{otherwise},
\end{cases}
\end{equation*}
where $V$ is again a $\kk$-central invertible $(S,S)$-bimodule concentrated in degree 0.

First, we note the following consequence of Proposition \ref{prop:extcalc}.
\begin{lemma}\label{lem:locgentleAStype} There exists $k$ such that $\grExt^{i}_A(S,A)=0$ for all $i\neq k$ if and only if either $\cQ$ is 2-regular or $A$ is one of $\kk\widetilde{A}_n$ and $(\kk\widetilde{A}_n)^\#$ for some $n$. 
Then, $\grExt^{k}_A(S,A)$ is an invertible $(S,S)$-bimodule concentrated in degree $-\ell$ where 
\begin{equation*}
    (k,\ell)=\begin{cases} 
    (0,-1) & \text{if }A=(\kk\tilde{A}_n)^\#,\\
    (1,1) & \text{if }A=\kk\tilde{A}_n,\\
    (1,0) & \text{if $\cQ$ is 2-regular.}
    \end{cases}
\end{equation*}
\end{lemma}

\begin{proof}
Suppose that there exists $k$ such that $\grExt^{i}_A(S,A)=0$ for all $i\neq k$.
By Corollary \ref{cor:depth},  $k=0$ if $A$ has a finite maximal path and $k=1$ otherwise. 
By Corollary \ref{cor:extcalc}, if $A^\#$ has a finite maximal path $\gamma$, then $\grExt^{\ell(\gamma)}_A(S,A)\neq 0$. Thus either  $A^\#$ has no finite maximal path (if $k=0$) or any finite maximal path of $A^\#$ has length 1 (if $k=1$). 

Suppose that $k=0$, so $A^\#$ has no finite maximal paths. By Lemma \ref{lem:invoutdeg}, this implies that either $A^\#=\kk\widetilde{A}_n$ or some vertex $u$ is the source of two arrows. By Corollary \ref{cor:extcalc}, the latter implies $\grExt^1_A(S(v),A)\neq 0$, so we must have $A^\#=\kk\widetilde{A}_n$. 

Now suppose that $k=1$, so $A$ has no finite maximal paths. Assume for contradiction that $A$ is not $\kk\widetilde{A}_n$ and $\abs{\cQ_1}<2\abs{\cQ_0}$. By Lemma \ref{lem:invoutdeg}, there must exist
vertices $u$ and $u'$ such that $\indeg(u)=\outdeg(u)=1$, $\indeg(u')=\outdeg(u')=2$ and there is some arrow $\alpha$ from $u$ to $u'$. 
As $t(\alpha)=u'$ and $u'$ is the source of two arrows, there must be some $\beta$ such that $\alpha\beta\in\cI$. As $\indeg(u)=1$, so $u$ is the target of exactly one arrow, $\beta'$, and $\beta'\alpha\notin\cI$. However, this means that $\gamma_\alpha^\#$ is finite and length at least two,
which is a contradiction. 
Thus, if $k$ exists, then $\cQ$ is 2-regular or $A$ is one of $\kk\widetilde{A}_n$ and $(\kk\widetilde{A}_n)^\#$.

Suppose $A=(\kk\widetilde{A}_n)^\#$. Using the convention of Example \ref{ex:Antilde}, $A=\kk\widetilde{A}_n/\langle\alpha_k\alpha_{k+1}\mid k=0,\ldots, n\rangle$. All maximal forbidden paths of $A$ are infinite and no vertex is the source of two arrows, so $\grExt^i_A(S,A)=0$ for $i\neq 0$ and 
\begin{equation*}
    \grExt^0_A(S,A)\cong \grHom_A(S,A)=\setst{\phi_{x}}{ x\in A_+},
\end{equation*}
where $\phi_{x}(y)\defeq xy,$  is concentrated in degree 1. Moreover, $\grHom_A(S,A)$ is an $(S,S)$-bimodule, with basis $\{\phi_{\alpha_i}\}$ and action given by $\epsilon_i\cdot \phi_{\alpha_j}=\delta_{i,j}\phi_{\alpha_j}$ and 
$\phi_{\alpha_j}\cdot \epsilon_i=\delta_{i,j+1}\phi_{\alpha_j}$. Let $V=\kk\{v_j\}_{1\leq j\leq n}$ be the $(S,S)$-bimodule given by $\epsilon_i\cdot v_j =\delta_{i,j+1} v_j$, $v_j\cdot \epsilon_i=\delta_{i,j}v_j$. Then
\begin{equation*}
   V \otimes_S \grHom_A(S,A)\cong\grHom_A(S,A)\otimes_S V\cong S,
\end{equation*}
so $\grExt^0_A(S,A)\cong \grHom_A(S,A)$ is an invertible $(S,S)$-bimodule, concentrated in degree 1.

If $A=\kk\widetilde{A}_n$, then $A$ is AS regular of dimension 1 \cite[Proposition 6.6]{RR19} and hence AS Gorenstein. To find the Gorenstein parameter $\ell$, we compute $\grExt^1_A(S(k),A)$.

By Proposition \ref{prop:extcalc}, 
\begin{equation*}
    \grExt^1_A(S(k),A)\cong \grHom_A(\epsilon_{k+1}A[-1],A)/\{\phi_x\}_{x\in A\alpha_k}\cong S(k+1)[1],
\end{equation*} 
as left $A$-modules, hence $\ell = 1$.

Lastly, if $\cQ$ is 2-regular, then neither $A$ nor $A^\#$ has a finite maximal path and, for each vertex $v$,  $\outdeg(v)=2$ and thus $\grExt^1_A(S(v),A)\cong \, _AS(v)$ as left $A$-modules by Proposition \ref{prop:extcalc}. 
The projective resolution of $S(v)$ is a complex of $(S,A)$-bimodules where the projective modules have the 
left $S$-module structure where $\epsilon_v$ acts as the identity and $\epsilon_{v'}$ acts as zero for any $v'\neq v$, so 
$\grExt^1_A(S(v),A)\cong S(v)$ as right $S$-modules as well. It follows that 
\begin{equation*}
    \grExt^1_A(S,A)\cong\bigoplus_{v\in \cQ_0}\grExt^1_A(S(v),A)\cong S.
\end{equation*}
\end{proof}

\begin{proposition}\label{prop:locgentleAScond}
    A (locally) gentle algebra is $A\neq \kk$ AS Gorenstein if and only if either $\cQ$ is 2-regular or $\cQ=\widetilde{A}_n$ and $A$ is  $\kk\widetilde{A}_n$ or $(\kk\widetilde{A}_n)^\#$. The only AS regular (locally) gentle algebra is   $\kk\widetilde{A}_n$.
\end{proposition} 

\begin{proof}
In Lemma \ref{lem:locgentleAStype}, we saw that no other (locally) gentle algebra can be AS Gorenstein or AS regular, as there exist $i\neq j$ such that $\grExt^i_A(S,A)$ and $\grExt^j_A(S,A)$ are both nonzero. By Theorem \ref{thm:injdim}, $A$ has (left and right) injective dimension 0 if $A=(\kk\widetilde{A}_n)^\#$ and 1 if $A=\kk\widetilde{A}_n$ or $\cQ$ is 2-regular. Thus, these three types of (locally) gentle algebras are AS Gorenstein, with Gorenstein parameter 1 if the quiver is $\widetilde{A}_n$ and 0 if $\cQ$ is 2-regular. However, if $\cQ$ is 2-regular or $A=(\kk\widetilde{A}_n)^\#$, then $A^\#$ has an infinite maximal path and thus $\gldim(A)=\infty$, so $A$ cannot be AS regular.
\end{proof}

Therefore, it follows from \cite[Theorem 6.3]{RR19} that the only twisted Calabi-Yau (locally) gentle algebra is $\kk\widetilde{A}_n$.

\subsection{Cohen-Macaulay}\label{subsec:CM}

There are multiple equivalent conditions for a commutative ring to be Cohen-Macaulay which are not necessarily equivalent in the non-commutative case. Thus, there are multiple potential definitions for a non-commutative ring to be Cohen-Macaulay. We consider definitions based on the following conditions.
\begin{rmk}[Cohen-Macaulay for a Commutative Ring]
Let $R$ be a commutative local Noetherian ring. Then $R$ is Cohen-Macaulay if the Krull dimension of $R$ is equal to its depth. If $R$ is affine, then $R$ is Cohen-Macaulay if the localization of $R$ at every maximal ideal is Cohen-Macaulay. If $R$ is connected graded, then it is well known that $R$ is Cohen-Macaulay if and only if there is a graded polynomial subalgebra $P$ of $R$ such that $R$ is a finitely generated, free $P$-module.
\end{rmk}

\begin{defn}
Let $\gA$ be a Noetherian graded algebra. Then $\gA$ satisfies
\begin{itemize}
    \item  \emph{(CM1)} if
    $\GKdim(A)=\dpth(\gA)$.
    \item  \emph{(CM2)} if there exist commuting, homogeneous, algebraically independent elements $x_1,x_2,\dotsc,x_d\in \gA$ so that $\gA$ is a finitely generated, free $\kk[x_1,\ldots, x_d]$-module and $d=\GKdim(\gA)$.
\end{itemize}
\end{defn}

In Theorem \ref{thm:CM}, we will show that (CM1) and (CM2) are equivalent for a (locally) gentle algebra, motivating a definition of Cohen-Macaulay for (locally) gentle algebras.

\begin{lemma}\label{lem:CMdef2} 
Let $A$ be a locally gentle algebra. 
Then $A$ satisfies (CM2) if and only if $A$ has no maximal paths of finite length. In this case, $A$ is a finitely generated free $\kk[x]$-module, where
$x\defeq \sum\limits_{\alpha\in\cQ_1} \alpha$ is
the sum of all arrows in the quiver $\cQ$.
\end{lemma}

\begin{proof} ($\Rightarrow$)
Since $\GKdim(A)=1$, $A$ satisfies (CM2) if and only if there is an element $x\in A$ such that $A$ is a finitely generated free $\kk[x]$-module. Let $x\defeq \sum\limits_{\alpha\in\cQ_1} \alpha$.
Observe that $x^n$ is the sum of paths of length $n$ in $A$, so $\epsilon_vx^n$ is the sum of the paths in $A$ with length $n$ and source $v$. Similarly, $\alpha x^n$ is the unique path of length $n+1$ which starts with arrow $\alpha$ if any such path exists, and zero otherwise. In particular, $\{\alpha x^n\}$ is linearly independent over $\kk$ if $\gamma_\alpha$ is infinite.

As $A$ is spanned by its paths, $A$ is generated as a $\kk[x]$-module by $\cQ_0\cup\cQ_1$. For each vertex $v\in\cQ_0$, choose an arrow $\alpha_v$ with source $v$, if such an arrow exists, and let $\widetilde{\cQ}_1$ denote the remaining arrows, $ \cQ_1\setminus\{\alpha_v\}_{v\in\cQ_0}$. 
For any $n\geq 0$, 
\begin{equation*}
    \alpha_vx^n=\epsilon_v x^{n+1}-\sum_{\substack{\alpha\in\widetilde{\cQ}_1\\ s(\alpha)=v}} \alpha x^n,\end{equation*} 
so $\cQ_0\cup\widetilde{\cQ}_1$ is also a generating set for $A$ as a $\kk[x]$-module.   

 Assume that $A$ has no finite maximal paths. We will show that  $\cQ_0\cup\widetilde{\cQ}_1$ is a $\kk[x]$-basis for $A$ and thus $A$ is a free $\kk[x]$-module. 
 Suppose that there are polynomials $p_y$ with coefficients in $\kk$ such that
 \begin{equation*}
     \sum_{y\in \cQ_0\cup\widetilde{\cQ}_1} yp_y(x)=0.
 \end{equation*}
  Fix $v\in\cQ_0$. As $A$ has no finite maximal paths,  the arrow $\alpha_v$ must exist. Moreover, as $\gamma_{\alpha_v}$ is infinite, there exists an arrow $\alpha'_v$ such that $\alpha'_v\alpha_v\notin\cI$. Since $\{\alpha'_v x^n\}$ is linearly independent over $\kk$ and
  \begin{equation*}
     0=\alpha'_v\sum_{y\in \cQ_0\cup\widetilde{\cQ}_1} yp_y(x)=\alpha'_v p_{\epsilon_v}(x),
 \end{equation*}
we must have $p_{\epsilon_v}=0$. By choice of $\widetilde{\cQ}_1$, there is at most one arrow $\alpha\in \widetilde{\cQ}_1$ with source $v$ and thus
 \begin{equation*}
     0=\epsilon_v\sum_{y\in \cQ_0\cup\widetilde{\cQ}_1} yp_y(x)=\epsilon_vp_{\epsilon_v}(x)+\alpha p_{\alpha}(x)=\alpha p_{\alpha}(x)
 \end{equation*}
 so $p_{\alpha}=0$ as well. Thus $p_y=0$ for each $y\in \cQ_0\cup\widetilde{\cQ}_1$ and hence $\cQ_0\cup\widetilde{\cQ}_1$ is a $\kk[x]$-basis for $A$.

($\Leftarrow$) Assume for contradiction that $A$ has at least one finite maximal path, $\gamma$, and that there exists an element $a\in A_+$ such that $A$ is a finitely generated free $\kk[a]$-module.  Then, there exist elements $y_1,\ldots, y_n$ which form a $\kk[a]$-basis for $A$ and polynomials $p_i$ such that $\gamma = \sum_{i=1}^n y_ip_i(a)$. Observe that 
\begin{equation*}
    0=\gamma a = \sum_{i=1}^n y_iap_i(a).
\end{equation*}
As the $y_i$ form a $\kk[a]$-basis for $A$, it must be that $ap_i(a)=0$ for each $i$. Since $\gamma$ is nonzero, $p_i(a)\neq 0$ for some $i$. Thus, as $ap_i(a)=0$ and $p_i(a)\neq 0$, we see that $\kk[a]$ is finite dimensional as a $\kk$-algebra. However, this is a contradiction to our assumption that $A$ is infinite dimensional over $\kk$ and finite dimensional over $\kk[a]$. Therefore, no such $x$ exists if $A$ contains any finite maximal paths.
\end{proof}

\begin{theorem}\label{thm:CM}
For a (locally) gentle algebra $A$, the following are equivalent.
\begin{enumerate}
    \item  $A$ satisfies (CM1).
    \item  $A$ satisfies (CM2).
    \item  $A$ contains no finite maximal paths or $A$ contains no infinite maximal paths.
\end{enumerate}
\end{theorem}

\begin{proof}
Note that the (locally) gentle algebra $A$ has no infinite maximal paths if and only if $A$ is gentle. If $A$ is gentle, $\GKdim(A)=\dpth(A)=0$ by Corollary \ref{cor:depth}, so $A$ satisfies (CM1). Since any finite dimensional $\kk$-algebra satisfies (CM2), every gentle algebra satisfies (CM2). 

Suppose $A$ is locally gentle. Then, $\GKdim(A)=1$ and, by Corollary \ref{cor:depth}, the depth of $A$ is 1 if $A$ has no finite maximal paths, so $A$ satisfies (CM1) if and only if $A$ has no finite maximal paths.  The result follows from  Lemma \ref{lem:CMdef2}.
\end{proof}

\begin{defn}
If a (locally) gentle algebra $A$ satisfies the equivalent conditions of Theorem \ref{thm:CM}, then we say that $A$ is \emph{Cohen-Macaulay}.
\end{defn}

If both $A$ and $A^\#$ are Cohen-Macaulay and locally gentle, then $\cQ$ must be 2-regular. 
If $A$ is Cohen-Macaulay and infinite dimensional but $\abs{\cQ_1}<2\abs{\cQ_0}$, then there is some vertex $v$ with $\indeg(v)=\outdeg(v)=1$ by Lemma \ref{lem:invoutdeg}. Let $\alpha$ and $\beta$ be the arrows with target and source $v$, respectively. As $A$ has no finite maximal paths, $\alpha\beta\notin\cI$ and thus $\gamma^\#_r(\alpha)=\alpha$, so $\gamma_\alpha^\#$ is finite. Hence either $A^\#$ is finite dimensional or $A^\#$ is not Cohen-Macaulay.

For commutative Cohen-Macaulay algebras, there is a strong connection between the Hilbert series and being Gorenstein. 
\begin{lemma}[Stanley's Theorem \cite{Stanley}] 
Suppose $R$ is a connected graded domain which is  Cohen-Macaulay. Then, $R$ is Gorenstein if and only if its Hilbert series satisfies
$h_R(t^{-1})=\pm t^{n}h_R(t)$ for some $n\in \ZZ$.
\end{lemma}
While all locally gentle algebras are Gorenstein (Theorem \ref{thm:injdim}), most do not satisfy this condition on their Hilbert series. Instead, we have the following analogue for locally gentle algebras.

\begin{theorem}\label{thm:Stanleylocgentlealgs}
Suppose that $A=\kk\cQ/\cI$ is locally gentle and Cohen-Macaulay. The following are equivalent.
\begin{itemize}
    \item[(i)] $A$ is AS Gorenstein.
    \item[(ii)] $h_A(t^{-1})=\pm t^{k}h_A(t)$ for some  $k\in\ZZ$.
    \item[(iii)] $\cQ$ is $\widetilde{A}_n$ 
    or $\cQ$ is 2-regular.
\end{itemize}
\end{theorem}

In particular, if $\cQ$ is $\widetilde{A}_n$, then $A=\kk\widetilde{A}_n$ and $h_A(t) = \frac{n+1}{1-t}=-t^{-1}h_A(t^{-1})$. If $\cQ$ is 2-regular, then  $h_A(t) = |\cQ_0|\frac{1+t}{1-t}=-h_A(t^{-1})$. In both cases, $k$ is the Gorenstein parameter of $A$ (1 and 0, respectively).

\begin{proof}
We saw that (iii) was equivalent to (i) in Lemma \ref{lem:h_Acond} and to (ii) in Proposition \ref{prop:locgentleAScond}. 
\end{proof}

These conditions are satisfied exactly when $A$ is semiprime (see Corollary \ref{cor:semiprime}). The condition that $A$ must be infinite dimensional is necessary, as we shall see in the following example.

\begin{example}\label{ex: Stanley fails}
Let $A$ be the path algebra of the Kronecker quiver, which is an orientation of  $C_2$. That is, $A$ is the path algebra of the following quiver.
\begin{center}
    \begin{tikzpicture}[scale=.7]
    \vtx{0,0}{1}
    \vtx{2,0}{2}
    \draw[->]  (.2,-.1) to node[midway,below]{$b$} (1.8,-.1) ;
    \draw[->] (.2,.1)  to node[midway,above]{$a$}  (1.8,.1);
    \end{tikzpicture}
\end{center}
Then $A=\kk\epsilon_1+\kk\epsilon_2+\kk a+\kk b$ is gentle and thus Cohen-Macaulay. This algebra is not AS Gorenstein since $\dpth(A)=0\neq 1=\injdim(A)$. In particular, $\grExt^0_A(S,A)\cong \grExt^0_A(S(2),A)\cong A\epsilon_2$ and
\begin{equation*}
    \grExt^1_A(S,A)\cong \grExt^1_A(S(1),A)\cong \left((A\epsilon_2\oplus A\epsilon_2)/\kk(a,b)\right)[1]\cong I(2),
\end{equation*}
so $A$ is not AS Gorenstein. However, $h_A(t)=2(t+1)=th_A(t^{-1})$.
\end{example}

The only AS Gorenstein gentle algebras are the algebras $A=(\kk\widetilde{A}_n)^\#$, which have Hilbert series $h_A(t)=(n+1)(1+t)$.

\subsection*{Acknowledgments} The authors would like to thank Jason Gaddis, Ken Goodearl, and Birge Huisgen-Zimmmerman for helpful discussions.
A. Oswald was supported by a grant from the Simons Foundation Targeted Grant (917524) to the Pacific Institute for the Mathematical Sciences. 
J.J. Zhang was partially supported by the US National Science Foundation (grant Nos. DMS-2302087 and  DMS-2001015). 

\printbibliography

\end{document}